\documentclass[11pt,reqno]{amsproc}

\title{The Benjamin--Feir stability diagram for Stokes waves with constant vorticity in finite depth}

\author{Bowen Chen}
\address{Department of Mathematics, University of Maryland, College Park, MD 20742}
\email{bchen0@umd.edu}

\author{Huy Q. Nguyen}
\address{Department of Mathematics, University of Maryland, College Park, MD 20742}
\email{hnguye90@umd.edu}



\usepackage{amsmath}
\usepackage{amsthm}
\usepackage{amssymb}
\usepackage{booktabs}
\usepackage{mathrsfs}
\usepackage[margin=1in]{geometry}
\usepackage{times}
\usepackage{setspace}
\usepackage{graphicx}
\usepackage{listings}
\usepackage{float}
\usepackage{tabularx}
\usepackage{stmaryrd}
\usepackage{subfigure}
\usepackage{listings}
\usepackage{verbatim} 
\usepackage{multirow}
\usepackage[table]{xcolor}
\definecolor{g}{rgb}{0,0.5,0}%
\usepackage{hyperref}

\usepackage{tikz}
\usetikzlibrary{decorations.pathmorphing}
\usepackage{mathtools}
\usepackage{bigstrut}
\usepackage{esint}
\usepackage{extarrows}
\usepackage{newtxtt}
\usepackage{times}
\usepackage{color}
\usepackage{adjustbox}

\theoremstyle{plain}
\newtheorem{theo}{Theorem}[section]
\newtheorem{prop}[theo]{Proposition}
\newtheorem{lemm}[theo]{Lemma}

\theoremstyle{definition}
\newtheorem{rema}[theo]{Remark} 

\numberwithin{equation}{section}
\newcommand{\bq}{\begin{equation}}
\newcommand{\eq}{\end{equation}}

\usepackage{tcolorbox}
\tcbset{
  colframe=black!75,
  colback=gray!10,
  coltitle=black,
  fonttitle=\bfseries,
  sharp corners,
  boxrule=0.5pt,
}
\newtcolorbox{mmain}[1][]{%
  colback=white,
  colframe=red!20,
  fonttitle=\bfseries,
  title={Input:},
  sharp corners,
  boxrule=0.5pt,
  #1
}

\newtcolorbox{mmaout}[1][]{%
  colback=gray!5,
  colframe=cyan!20,
  fonttitle=\bfseries,
  title={Output:},
  sharp corners,
  boxrule=0.5pt,
  #1
}

\newtcolorbox{mmacode}[1][]{%
  colback=white,
  colframe=green!20,
  fonttitle=\bfseries,
  title={Code:},
  sharp corners,
  boxrule=0.5pt,
  #1
}

\DeclareMathOperator{\Span}{span}

\renewcommand{\th}{\mathtt{t}_h}

\newcommand{\R}{\mathbb{R}}

\def\L{\mathcal{L}}
\def\Ls{\mathscr{L}}
\def\sign{\mathrm{sign}}
\def\Z{\mathbb{Z}}
\def\C{\mathbb{C}}
\def\pinv{\p_x^{-1}}

\def\Lmat{\mathtt{L}}
\def\res{\operatorname{Res}_{\lambda=0}}
\def\anaR#1{\R\{#1\}}
\def\anaC#1{\C\{#1\}}
\newcommand{\vet}[2]{\begin{bmatrix}#1 \\ #2 \end{bmatrix}}
\def\B{\mathcal{B}}
\def\coeffM{\mathtt{c}_\mu}
\def\coeffE{\mathtt{c}_\eps}
\def\Os{\Omega_\sharp}
\def\Oscoeff#1{\Os^{(#1)}}

\def\fe{\mathtt{f}_\eps}

\def\hcrit{h_{\text{crit}}}
\def\dcrit{d_{\text{crit}}}
\def\gammacrit{\gamma_{\text{crit}}}

\def\eps{\varepsilon}
\def\na{\nabla}

\def\wh{\widehat}

\def\mez{\frac{1}{2}}

\def\Rr{\mathbb{R}}
  
\def\Nn{\mathbb{N}}
\def\Zz{\mathbb{Z}}
\def\Cc{\mathbb{C}}
\def\cR{\mathcal{R}}

\def\cK{\mathcal{K}}

\def\cL{\mathcal{L}}
\def\cR{\mathcal{R}}

\def\cV{\mathcal{V}}
\def\ld{\lambda}

\def\p{\partial}

\def\na{\nabla}

\def\ol{\overline}
\def\T{\mathbb{T}}

\def\wt{\widetilde}

\def\Om{\Omega}
\def\ld{\lambda}

\def\om{\omega}

\newcommand{\rem}[1]{O(\eps^{#1})}

\begin{document}
\maketitle
\begin{abstract}
We determine the Benjamin--Feir stability diagram for small-amplitude Stokes waves with constant vorticity in finite depth. For the full water wave system, we derive an explicit instability index and rigorously classify its sign throughout the depth-vorticity plane. With gravitational acceleration and carrier wavenumber normalized to one, we identify a unique minimum depth $h_{\text{crit}}\approx1.36271$, below which Benjamin--Feir instability is absent for every vorticity. Above this depth, two vorticity thresholds delimit the unstable regime,  outside of which the spectrum near the origin remains purely imaginary for all small Floquet exponents. In particular, at each fixed finite depth, constant vorticity of sufficiently large magnitude suppresses Benjamin--Feir instability. Moreover, any constant vorticity  $\gamma\ge  2/\sqrt{3}$ suppresses Benjamin--Feir instability at every finite depth.  The proof combines Hamiltonian spectral perturbation theory, symbolic computations, and formal verification in Rocq. Away from the transition set, our results rigorously confirm the stability criterion predicted by the nonlinear Schr\"odinger approximation of Thomas, Kharif, and Manna.
\end{abstract}
\tableofcontents
\section{Introduction}
This work concerns the stability of the classical small-amplitude Stokes waves \cite{Stokes1847, Nekrasov, Levi_Civita, Struik1926}. The surprising discovery of  Benjamin and Feir \cite{BenjaminFeir1967} in 1967 showed that {\it irrotational} Stokes waves are unstable with respect to  longitudinal, long-wave perturbations. This Benjamin--Feir/modulational instability has been proven to be ubiquitous in nonlinear wave phenomena \cite{ZakharovOstrovsky2009}. In the context of classical Stokes waves,  the first rigorous proof of Benjamin--Feir  instability was found for the finite-depth case by Bridges and Mielke \cite{finite_depth_BF_instability_1995} in 1995. The infinite depth case was recently resolved by Nguyen and Strauss \cite{Huy21}. In both \cite{finite_depth_BF_instability_1995} and \cite{Huy21}, unstable eigenvalues near $0\in \Cc$ were rigorously constructed. On the other hand, numerical computations have revealed that the Benjamin--Feir spectrum forms a figure-$8$ curve center at $0\in \Cc$. This full Benjamin--Feir spectrum was rigorously justified by Berti, Maspero, and Ventura in \cite{Berti_Full, Berti_finite_depth, critical_depth_unstable}, for both infinite and finite depth cases. Alternative proofs of the Benjamin--Feir instability were obtained in \cite{Yang} and \cite{Hur}. We also refer to  \cite{Hur, BertiMasperoVentura2025, Berti_isola} for the high-frequency instability  and  \cite{Huy_transverse, Huy_transverse_finite_depth, Transverse-RSY, Berti-McLean} for the transverse instability of irrotational gravity Stokes waves. Regarding irrotational capillary-gravity Stokes waves, we refer to \cite{HurYang-capillary, HsiaoMaspero2025, HsiaoWang2026} for Benjamin--Feir and high-frequency instabilities, and \cite{Truong_transverse} for transverse instabilities.

In this work, we investigate the stabilizing/destabilizing effect of constant vorticity on the Benjamin--Feir (modulational) instability of periodic traveling waves. The constant vorticity $\gamma$ is generated by a background current with a linear velocity profile, ${\bf u_s}=-\gamma y\mathbf e_x$. Such shear currents occur in coastal and tidal flows and in wind-driven surface layers \cite{Lentz1992, StaceyMonismithBurau1999}. To understand the effect of constant vorticity on periodic wave-trains,  Thomas, Kharif and Manna \cite{ThomasKharifManna2012} formally derived a nonlinear Schr\"odinger (NLS) equation governing the slowly varying wave envelope, accounting for the coupling between the background vorticity and the mean flow induced by modulation. Their NLS has the form 
\bq\label{NLS}
ia_\tau+La_{\xi\xi}=M|a|^2a,
\eq
where the constants $L\in \Rr$ and $M\in \Rr$  depend on the constant vorticity $\gamma\in \Rr$ and the fluid's depth $h\in (0, \infty]$. Considering the  stability of  plane-wave solutions $a_*= a_0\exp(-iMa_0^2\tau)$, $a_0\in \Rr$,  the authors found that $a_*$ is stable if \eqref{NLS} $LM>0$ (defocusing) and unstable  if  $LM<0$ (focusing). In the context of the formally derived NLS, the instability index $LM$ yields a stability diagram in the depth-vorticity plane.  In particular, in infinite depth, the plane-wave solution $a_*$ undergoes a transition from modulational instability to stability as $\gamma$ increases through the {\it critical value} 
\bq\label{gammacri}
\gamma_*=2\sqrt\frac{gk}{3},
\eq
 where $g>0$ is the gravitational acceleration and $k$ is the carrier wavenumber. In our earlier work \cite{ChenNguyen}, we rigorously esbtablished the transition threshold \eqref{gammacri} for the full water wave system with constant vorticity in infinite depth. 
 
The present paper concerns the finite depth case. Our main result is the rigorous justification of the instability index $LM$ of Thomas, Kharif and Manna \cite{ThomasKharifManna2012} for the full water wave system with constant vorticity in any finite depth $h\in (0, \infty)$. In order to state the result, we fix a constant vorticity $\gamma$ and a corresponding  Stokes wave of small amplitude $\eps$ and wave number $k=1$ (without loss of generality). We  perturb the Stokes wave by longitudinal, long-wave perturbations  of Bloch form \(e^{i\mu x}f(x)\), where $f$ is $2\pi$-periodic and $\mu$ is small.  After linearizing the full water wave system around the Stokes wave and performing a number of transformations, we arrive at a linear Hamiltonian operator $\L_{\mu,\eps}: H^1(\T)^2\to L^2(\T)^2$, $\T:=\Rr/2\pi \Zz$.  We refer to  \eqref{def:analytic_r} for the notation $r(\mu^{m_1}\eps^{n_1},\dots,\mu^{m_p}\eps^{n_p})$ for analytic remainder functions. After normalizing  $g=1$, our main result is the following complete description of the spectrum of $\L_{\mu,\eps}$ near the origin in the complex plane.
\begin{theo}\label{theo:main1}
Let $\gamma \in \Rr$, $h\in (0, \infty)$, and denote
\bq\label{def:c0:0}
\th=\tanh(h),\quad c_0=\mez\left(\gamma \th+\sqrt{\gamma^2\th^2+4\th}\right).
\eq
 There exist $\eps_0$, $\mu_0>0$, and an open neighborhood $U$ of $0\in \Cc$ such that the following holds. For all $\eps\in (-\eps_0, \eps_0)$, $\mu\in (0, \mu_0)$, the spectrum of the linearized operator $\L_{\mu, \eps}$ in $U$ consists of four eigenvalues $\ld^*_\pm$ and $\wt{\ld}_\pm$, where $\ld^*_\pm\in i\Rr$ and  

\bq
\wt{\ld}_\pm=\frac{i\mu(c_0^3h\th^2-c_0^3h+c_0^3\th)}
             {\th(c_0^2+\th)}
      +i\mu r(\mu^2,\mu\eps,\eps^2)\pm \mez\mu \sqrt{\Delta},
\eq
\bq
\Delta:=\coeffM\mu^2+\coeffE\eps^2
      +\wt{r}(\mu^3,\mu^2\eps,\mu\eps^2,\eps^3).
      \eq
Here $\coeffM=\coeffM(h, \gamma)<0$ is given by \eqref{def:cmu}, and $\coeffE=\coeffE(h, \gamma)$ is given by \eqref{def:ceps} with $\sign(\coeffE)=-\sign(\Omega)$, $\Omega$ depending only on $h$ and $\gamma$. 

(i) (Absence of Benjamin--Feir instability) If $(h, \gamma)$ satisfies $\Omega(h, \gamma)>0$, then the Stokes wave with amplitude $\eps$ does not exhibit Benjamin--Feir instability.  

(ii) (Benjamin--Feir instability)  Assume that $(h, \gamma)$ satisfies $\Omega(h, \gamma)<0$.  Then there exists a piecewise analytic function $\underline{\mu}:(-\eps_0,\eps_0)\to [0,\mu_0)$ depending on $(h, \gamma)$, of the form
    $$
\underline{\mu}(\eps)=
    \begin{cases}
        \sqrt{-\frac{\coeffE}{\coeffM}}\eps+r_+(\eps^2) &\quad\text{if~} \eps\ge 0,\\
        - \sqrt{-\frac{\coeffE}{\coeffM}}\eps+r_-(\eps^2) & \quad\text{if~} \eps< 0,
    \end{cases}
    $$ 
    such that for any $\eps\in(-\eps_0,\eps_0)$, we have 
    \bq
    \Delta>0~\text{if}~0<\mu<\underline{\mu}(\eps),\quad\Delta<0~\text{if}~\underline{\mu}(\eps)<\mu<\mu_0.
    \eq
\end{theo}
We have verified using Mathematica that  the index  $LM$ in \cite{ThomasKharifManna2012} can be written as
\bq
    LM=-\frac{1}{8}\coth^2(h)\,\coeffE.
\eq
Therefore, the NLS instability condition $LM<0$ is equivalent to our  instability condition $\Omega<0$ since $\sign(\coeffE)=-\sign(\Omega)$. In the unstable regime $\Omega(h, \gamma)<0$, as $\mu$ varies in $(0, \underline{\mu}(\eps))$, the eigenvalues $\wt{\ld}_\pm$ trace the upper half of a  figure-$8$ curve center at $0\in \Cc$; the lower half comes from $\sigma(\cL_{-\mu, \eps})=\overline{\sigma(\cL_{-\mu, \eps})}$.  

The index $\Omega(h, \gamma)$ written in terms of $(h, c_0^2)$ is the sextic polynomial $\Omega_\sharp(\gamma, d)$ of $d=c_0^2$, given by \eqref{eq:Omega_explicit}-\eqref{form:Omegaj}. We use formal proofs verified in Rocq to establish the global sign classification of  $\Omega$. This yields the following quantitative description of the modulationally stable and unstable regions in the $(h, \gamma)$-plane. 
\begin{theo}[The sign of $\Omega$]\label{theo:main2}
The following statements hold.
\begin{enumerate}
    \item
    There is a unique point
    $(\hcrit,\dcrit)\in(0,\infty)^2$
    satisfying $  \Omega_\sharp(\hcrit,\dcrit)
        =\partial_d\Omega_\sharp(\hcrit,\dcrit)=0$. 
    Moreover,
    \[
        1.36271<\hcrit<1.36272,
        \qquad
        0.86010<\dcrit<0.86011.
    \]
  Denote 
    \[
        \gammacrit
        :=\frac{\sqrt{\dcrit}}{\tanh(\hcrit)}
          -\frac{1}{\sqrt{\dcrit}}
        \in(-0.02081,-0.02079).
    \]
     \item
    If $0<h<\hcrit$, then $ \Omega(h,\gamma)>0$ for all $\gamma\in\mathbb R$. At the critical depth, $\Omega(\hcrit,\gamma)\ge 0$ for all $\gamma\in\Rr$, with equality if and only if $\gamma=\gammacrit$.
    \item Consider $h>\hcrit$. The equation
    $\Omega(h, \cdot)=0$ has exactly two real roots,
    denoted by $\gamma_l(h)<\gamma_r(h)$, which are smooth functions of $h$. The sign of $\Omega$ is given by
    \[
    \begin{cases}
        \Omega(h,\gamma)>0
        & \text{if }\gamma<\gamma_l(h)
          \text{ or }\gamma>\gamma_r(h),\\
        \Omega(h,\gamma)=0
        & \text{if }\gamma=\gamma_l(h)
          \text{ or }\gamma=\gamma_r(h),\\
        \Omega(h,\gamma)<0
        & \text{if }\gamma_l(h)<\gamma<\gamma_r(h).
    \end{cases}
    \]
    \item For any $\gamma\ge \frac{2}{\sqrt3}$ and any finite depth $h$, $\Om(h,\gamma)>0$. 
\end{enumerate}
\end{theo}
We refer to Figure \ref{figure:Om} below for a plot, generated using Mathematica, of the zero level set $\Omega=0$ and the modulational stability and instability regions it separates.

We recall that the Whitham-Bejamin critical depth is $h_{\text{WB}}\approx 1.363$ \cite{finite_depth_BF_instability_1995, Berti_finite_depth, Hur}, below which all small-amplitude {\it irrotational} Stokes waves are modulationally stable. By part (2) of Theorem \ref{theo:main2}, there is slightly smaller critical depth $\hcrit\approx 1.3627<h_{\text{WB}}$ below which small-amplitude Stokes waves with arbitrary constant vorticity are modulationally stable. On the other hand, part  (3) of Theorem \ref{theo:main2} shows that for any finite depth above $\hcrit$, constant vorticity of sufficiently large magnitude suppresses Benjamin--Feir instability. Finally, any constant vorticity above  $\gamma_*=\frac{2}{\sqrt{3}}$ suppresses Benjamin--Feir instability at every finite depth. Combining this with the result in  \cite{ChenNguyen} for the infinite-depth case, we conclude that 
\begin{align*}
&\textit{Any constant vorticity $\gamma>\frac{2}{\sqrt{3}}$ suppresses Benjamin--Feir instability}\\
& \qquad \textit{of small-amplitude Stokes waves at any depth $h\in (0, \infty]$. }  
\end{align*}
The proof of  Theorem \ref{theo:main1} follows the strategy in our infinite-depth work \cite{ChenNguyen}. Let us explain the key steps and highlight differences. From \cite[Theorem 1.6]{ChenNguyen}, for any constant vorticity $\gamma$, we have a family of Stokes waves with small amplitude $\eps$ and the expansion of them in second order of $\eps$. We fix such a Stokes wave and  linearize the generalized Zakharov-Craig-Sulem  system \eqref{eq:main_with_t} around it.  Using the linear transformation due to Wahl\'en \cite{Wahlen_Hamiltonian}, we transform the linearized system to the canonical Hamiltonian form.  After using Alinhac's good-unknown and a conformal mapping, we introduce the small Floquet exponent $\mu$ and arrive at the linearized operator $\L_{\mu, \eps}$ \eqref{def:finalL}. In the infinite depth case \cite{ChenNguyen}, after similar transformations, we obtained the singular opertor $(\mu+D)^{-1}\sign(D+\mu)$. In the finite depth case, however, this term becomes $\frac{\tanh\big((h+\fe)(D+\mu)\big)}{D+\mu}$, which is analytic in $(\mu, \eps)$ near $(0, 0)$. When $\mu=0$, $\ld=0$ is an eigenvalue with algebraic multiplicity four and an explicit kernel $\cV_{0, \eps}$. Following \cite{Berti_Full}, we use Kato's similarity transformation to obtain a formula for the perturbed basis $\cV_{\mu, \eps}$ in terms of the unperturbed basis $\cV_{0, 0}$ and the spectral projection of $\cL_{\mu, \eps}$.  Consequently, eigenvalues near $0\in \Cc$  of $\L_{\mu, \eps}$ are precisely eigenvalues of the $4\times 4$ matrix $\Lmat_{\mu,\eps}$ representing  $\L_{\mu, \eps}$ on the $\mathcal{V}_{\mu, \eps}$ basis. We note that the corresponding representation matrix in \cite{ChenNguyen} is a $3\times 3$ matrix. Combining the Hamiltonian structure and the reversibility of $\L_{\mu, \eps}$ yields symmetries of entries of $\Lmat_{\mu,\eps}$. We expand $\Lmat_{\mu,\eps}$ to second order in $\mu$ and $\eps$. After that, we prove exact vanishing of certain entries in the $\mu^0$ and $\mu^1$ coefficients of  $\Lmat_{\mu,\eps}$. The vanishing is exact in the sense that it is valid at all orders of $\eps$. To this end, we utilize the aforementioned spectral property of $\cL_{0, \eps}$ and parity properties. With the above expansion and exact vanishing, we implement the direct method in \cite{ChenNguyen} to analyze the spectrum of  $\Lmat_{\mu,\eps}$ and obtain Theorem \ref{theo:main1}. 

After completing this work, we became aware of the independent works \cite{FanSun2026, HsiaoHurYang2026}, which establish closely related results. The present work was developed independently and is based on our infinite-depth paper \cite{ChenNguyen}. 
\subsection{Notation}\label{sec:notation}
  \leavevmode
\begin{itemize}
\item $\th=\tanh(h)$ for $h\in (0, \infty)$.
    \item $\T_L:=\R/(L\Z)$, $\T\equiv \T_{2\pi}$.
    \item The Fourier transform for an $L$-periodic function is defined by
    \begin{equation}
        \wh{f}^L(k):=\int_{\T_L} f(x)e^{-i\frac{2\pi}{L}kx}dx,\quad \forall k\in\Z,\quad \wh{f}(k) \equiv \wh{f}^{2\pi}(k),    \end{equation}
    so that
    \begin{equation}
        f(x)=\frac{1}{L}\sum_{k\in\Z}\wh f^L(k)e^{i\frac{2\pi}{L}kx},\quad f\in L^2(\T_L).
    \end{equation}
    Fourier multipliers are then defined by
\begin{equation}\label{def:FM}
    m\left(D_L\right) f(x)=\frac{1}{L} \sum_{k \in \mathbb{Z}} e^{i \frac{2 \pi}{L} k x} m\left(k \frac{2 \pi}{L}\right) \hat{f}^L(k),\quad m(D)\equiv m(D_{2\pi}).
\end{equation}
    \item For $\sigma,s\ge 0$ and $f: \T\to \Cc$,  we define the norm
        \begin{align}
        \|f\|_{H^{\sigma,s}}^2&:=\sum_{k\in\Z}(1+|k|^2)^{s}e^{2\sigma|k|}|\wh f(k)|^2
    \end{align}
   and the corresponding function space
        \begin{align}
            H^{\sigma,s}(\T):=\{f\in L^2(\T):~\|f\|_{H^{\sigma,s}}<\infty\}.
        \end{align}
    \item For $f,g \in \left(L^2(\T)\right)^2$, 
    $$(f, g):=\int_{\T} f(x) \cdot ~\overline{g}(x)dx,$$ 
    where the `bar' denotes complex conjugation  and the `dot' denotes the real dot product. 
    \item For an $L$-periodic function $f$, we define the projections onto the zero mode and nonzero modes by
\[
\Pi_0f=\fint f=\frac{1}{L}\int_0^Lfdx,\quad \Pi_0^\perp f=f-\Pi_0f.
\]
If $\Pi_0f=0$, we denote the unique mean-free primitive of $f$ by $\pinv f$, which is $L$-periodic.  We note the identities  
\begin{equation}\label{dinverse}
\p_x\pinv=\mathrm{Id},\quad \pinv\p_x=\Pi_0^\perp.
\end{equation}
\item For independent variables $(X_1,\dots,X_n)\in\R^n$ and a Banach space $B$, we denote the space of $B$-valued analytic functions near $0\in \R^n$ by $B\{X_1,\dots,X_n\}$. 
\begin{itemize}
    \item For a function $F(X_1,\dots,X_n)\in B\{X_1,\dots,X_n\}$, we denote by $[X_j^m]F$ the coefficient of $X_j^m$ in the expansion of $F$, that is, 
    \begin{equation}\label{nota:coefficient}
    [X_j^m]F:=\frac{\p_{X_j}^mF|_{X_j=0}}{m!}.
   \end{equation}
\item If $F(\mu,\eps)$ is an analytic function in a neighborhood of $(0, 0)$, by writing $F=O(\mu^{m_1}\eps^{n_1},\dots,\mu^{m_p}\eps^{n_p})$ we mean $F(\mu,\eps)=\sum_{j=1}^p \mu^{m_j}\eps^{n_j}G_j(\mu,\eps)$ for some analytic functions $G_j(\mu,\eps)$ near $(0,0)$.
Note that $O(\mu^{m_1}\eps^{n_1},\dots,\mu^{m_p}\eps^{n_p})$ defined here is not  the same as the standard notation $F=O(|\mu^{m_1}\eps^{n_1}|+\dots+|\mu^{m_p}\eps^{n_p}|)$, where the latter means 
\[
|F|\le C\sum_{k=1}^p |\mu^{m_k}\eps^{n_k}|.
\]
  For example,  $\mu \eps =O(\mu^2+\eps^2)$  but $\mu\eps \ne O(\mu^2,\eps^2)$.
\item By $r(\mu^{m_1}\eps^{n_1},\dots,\mu^{m_p}\eps^{n_p})\in\anaR{\mu,\eps}$ we mean any real analytic function satisfying 
   \begin{equation}
   \label{def:analytic_r}
       r(\mu^{m_1}\eps^{n_1},\dots,\mu^{m_p}\eps^{n_p})=O(\mu^{m_1}\eps^{n_1},\dots,\mu^{m_p}\eps^{n_p}).
   \end{equation}
\end{itemize}
   \item We denote by $even$ (respectively, $odd$) a $2\pi$-periodic even (respectively, odd) function  on $\T$. They may denote different functions in different contexts.
\end{itemize}
\subsection{Water waves with constant vorticity} 
We consider two-dimensional inviscid, incompressible fluids in the moving domain  $\Omega(t)=\{(x, y):x \in \mathbb{R}, -h< y<\eta(x, t)\} $ of finite depth $h\in (0, \infty)$, where $\{y=\eta(x,t)\}$ is the free surface. We consider a shear current   ${\bf u_s}=-\gamma y \mathbf e_x$ which generates constant vorticity $\gamma=\nabla^\perp\cdot {\bf u_s}$, where $\nabla^\perp=(-\p_y, \p_x)$ and $\mathbf e_x$ is the unit vector in the $x$-direction. Then, we consider fluid velocity of the form  
\begin{equation}
{\bf u}=
        -\gamma y \mathbf e_x+\na\Phi,
    \end{equation}
     where  $\Phi$ is a harmonic  potential, so that $\nabla^\perp\cdot{\bf u}=\gamma$. Denoting $\psi(x, t)=\Phi(x, \eta(x, t), t)$, the potential $\Phi$ can be recovered from $\eta$ and $\psi$ as the solution of the Laplace equation
     \begin{equation}
     \begin{cases}
     \Delta\Phi=0\quad\text{in~}\Omega,\\
     \Phi(x, \eta(x))=\psi(x),\quad \p_y\Phi(x, -h)=0.
     \end{cases}
     \end{equation}
     We denote the Dirichlet-Neumann operator $G_h(\eta)\equiv G(\eta)$ by 
     \begin{equation}
     G_h(\eta)\psi(x)=-\eta_x(x)\Phi_x(x, \eta(x))+\Phi_y(x, \eta(x)).
     \end{equation}
      Using $G_h$, the water wave problem can be recast as a coupled system of $(\eta, \psi)$:
\begin{equation}
\label{eq:main_with_t}
\begin{cases}
\eta_t=F_1(\eta,\psi,c):= c\partial_x\eta +G_h(\eta) \psi+\gamma \eta \eta_x \\
\psi_t=F_2(\eta,\psi,c,P):= c\partial_x\psi -g \eta-\frac{{\psi_x^2}}{2}+\frac{\big(\eta_x \psi_x+G_h(\eta) \psi\big)^2}{2\left(1+\eta_x^2\right)}+\gamma \eta \psi_x+\gamma \partial_x^{-1} G_h(\eta) \psi+P,
\end{cases}
\end{equation}
where $c\in \Rr$ is the speed of a horizontal moving frame,  $g>0$ is gravitational acceleration and $P\in \Rr$ is Bernoulli constant. \eqref{eq:main_with_t} is a generalization of the well-known Zakharov-Craig-Sulem formulation \cite{Zakharov, Craig_Sulem} for irrotational flows. We refer to \cite{ConstantinIvanovProdanov2008, Wahlen_Hamiltonian} for the derivation of \eqref{eq:main_with_t}. 

We normalize $g$ to one by the change of variables 
\begin{equation}
\label{eq:nondim}
    (\eta, \psi, c, P, t,  \gamma)\mapsto(\eta, \sqrt{g}\psi, \sqrt{g}c, gP, \frac{t}{\sqrt{g}}, \sqrt{g}\gamma),
\end{equation}
Any time-independent periodic solution of \eqref{eq:main_with_t} is called a Stokes wave with speed $c$. The existence of small-amplitude Stokes waves with constant vorticity has been proven in \cite{ConstantinVarvaruca2011, Berti_bifurcation_with_surface_tension}. The following theorem records both the existence and second-order expansion of these Stokes waves.
\begin{theo}\cite[Theorem 1.6]{ChenNguyen}
\label{thm:Stokes_without_g}
Let  $\sigma>0$ and $s\ge\frac72$ such that $s+\mez \in\Z$. For any constant vorticity $\gamma$ and any depth $h\in (0, \infty)$, there exist $\eps_{St}=\eps_{St}(\gamma,h)>0$ and a unique family of solutions $\big(\eta(\cdot; a,h), \psi(\cdot; a,h), c(a,h)\big)\in H^{\sigma,s}(\T)\times H^{\sigma,s}(\T)\times \Rr$ to the steady equation \eqref{eq:main_with_t} with $P=0$, parametrized by $a\in (-\eps_{St}, \eps_{St})$ such that 

(i) the mapping $a\mapsto \big(\eta(\cdot; a,h), \psi(\cdot; a,h), c(a,h)\big)\in H^{\sigma, s}(\T)\times H^{\sigma, s}(\T)\times \Rr$ is analytic, and

(ii) $\eta(\cdot; a,h)$ is even and $\psi(\cdot; a,h)$ is odd. 

Moreover,  the following expansions hold: 
\begin{equation}
    \label{expand:Stokes_without_g_all_depth}
        \begin{aligned}
c=c(a,h)&=c_0+a^2c_2+O(a^3),\\ 
    \eta(x)=\eta(x;a,h)&= 
        a \Big[\th \cos(x)\Big]+\frac{ a ^2}4\Big[(c_0^2-2\th+c_0^2\th ^2)\\
        &+\Big(\frac{1}{\th}+\frac{c_0^2+\th(1-2c_0^4)+c_0^2\th ^2}{c_0^4}\Big)\cos(2x)\Big]+O(a^3),
    \\
    \psi(x)=\psi(x;a,h)&= a \Big[ c_0\sin(x)\Big]+ a ^2\Big[\frac{c_0^4+c_0^2\th+(1-2c_0^4)\th^2+c_0^2\th^3+(1+c_0^4)\th^4}{8c_0^3\th^2}\sin(2x)\Big]+O(a^3),
\end{aligned}
\end{equation}
and 
\begin{equation}
    \begin{aligned}
        c_2=c_2(h)&=\frac{1}{8 c_0^5 \th^2 \left(c_0^2+\th\right)}\Bigg(c_0^8+2 c_0^6 \th+\left(11 c_0^8-7 c_0^4+1\right) \th^4+\left(3 c_0^4-5 c_0^8\right) \th^2\\
        &+\left(-2 c_0^{10}+4 c_0^6+2 c_0^2\right) \th^5+\left(-2 c_0^{10}-2 c_0^6+2 c_0^2\right) \th^3\Bigg),
    \end{aligned}
\end{equation}
 where $\th:=\tanh(h)$ and $c_0$ is the positive solution of 
\begin{equation}
    \label{eq:relation_c0_gamma_general}
\frac{c_0}{\th}-\frac{1}{c_0}=\gamma.
\end{equation}
\end{theo}
The general case $P\in \Rr$ can be reduced to the case $P=0$ by setting
\begin{equation}
    \label{eq:symm_P}
    \wt{\eta}=\eta-P,\quad \wt{c}=c+\gamma P, \quad \wt h=h+P
\end{equation}
and using the fact that $G_h\left(\widetilde{\eta}+P\right) \psi=G_{\wt h}(\widetilde{\eta}) \psi$.  Therefore, from the family of Stokes waves  \eqref{expand:Stokes_without_g_all_depth}, we obtain the following family of Stokes waves with general $P\in \Rr$:
\begin{equation}
    \label{eq:solution_with_P}
    \breve\eta(x;a,h,P)=\eta\Big(x;a,h+P\Big)+P,\quad 
    \breve\psi(x;a,h,P)=\psi\Big(x;a,h+P\Big),\quad 
   \breve c(a,h,P)=c\Big(a,h+P\Big)-\gamma P.
\end{equation}
\section{Linearization and transformations}
We fix  a Stokes wave $\left(\eta^*, \psi^*, c^*, P^*\right)$ as given in \eqref{expand:Stokes_without_g_all_depth} with  $a=\eps$ and $P^*=0$.  We consider $L$-periodic perturbations $(\ol{\eta}, \ol{\psi})$ of $(\eta^*, \psi^*)$, where $L=n_02\pi$ for some $n_0\in \Nn$. We linearize \eqref{eq:main_with_t} around $(\eta^*, \psi^*)$ and use Alinhac's  good unknown
\begin{equation}
    v_1=\ol{\eta}, \quad v_2=\ol{\psi}-B^* \ol{\eta},
\end{equation}
where 
\begin{equation}
\label{def:B_and_V}
    B^*=\frac{G(\eta^*) \psi^*+\partial_x \psi^* \partial_x \eta^*}{1+\left|\partial_x \eta^*\right|^2}, \quad V^*=\partial_x \psi^*-B^* \partial_x \eta^*. 
\end{equation}
From \cite[Section 2.1]{ChenNguyen}, we obtain 
\begin{align}
    \label{eq:v1}\p_tv_1&=\p_x\Big((c^*-V^*+\gamma\eta^*)v_1\Big)+G(\eta^*)v_2,\\
    \label{eq:v2}\p_tv_2&=-\Big(1+(V^*-c^*-\gamma\eta^*)\p_xB^*\Big)v_1+\gamma \Pi_0(V^* v_1)+(c^*-V^*+\gamma\eta^*)\p_xv_2+\gamma\p_x^{-1}G(\eta^*)v_2.
\end{align}
In matrix form,
\begin{equation}
    \label{eq:v_matrix}
    \begin{aligned}
    \p_t \vet{v_1}{v_2}&=\Ls_\eps\vet{v_1}{v_2}\\
    &:=
    \begin{bmatrix}
        \p_x\circ (c^*-V^*+\gamma\eta^*)& G(\eta^*)\\
        -\Big(1+(V^*-c^*-\gamma\eta^*)\p_xB^*\Big)+\gamma\Pi_0\circ V^* & (c^*-V^*+\gamma\eta^*)\p_x+\gamma\p_x^{-1}G(\eta^*)
    \end{bmatrix}
    \vet{v_1}{v_2}.
    \end{aligned}
\end{equation}
 As shown in \cite[Theorem 2.4]{ChenNguyen}, which remains valid in finite depth, we have
\begin{theo}
        \label{thm:eigenvectors_R}
 For all $\gamma \in \Rr$, there exists $\eps_0=\eps_0(\gamma)>0$ such that for all $0<|\eps|<\eps_0$,  $0$ is an eigenvalue of $\Ls_\eps$ with algebraic multiplicity $4$ and geometric multiplicity $2$ or $3$. Moreover, 
    \begin{align}
        R_1&=[0,1]^T,\\
        R_2&=\Big[\p_x\eta^*, \p_x\psi^*-B^*\p_x\eta^*\Big]^T
    \end{align}
    are eigenvectors, 
\begin{equation}
     R_3=\Big[\p_a\eta\vert_{a=\eps}, \big(\p_a\psi-B^*\p_a\eta\big)\vert_{a=\eps}\Big]^T
     \end{equation}
     can be an eigenvector or generalized eigenvector, depending on $(\eps, \gamma, h)$,  and 
               \begin{equation}
                R_4=\Big[\p_P\breve\eta\vert_{(a, P)=(\eps, 0)}, \big(\p_P\breve\psi-B^*\p_P\breve\eta\big)\vert_{(a, P)=(\eps, 0)}\Big]^T
\end{equation}
is a generalized eigenvector. Moreover,  
\bq\label{LR3R4}
\Ls_\eps R_3,\Ls_\eps R_4\in\Span\{R_1,R_2\}.
\eq
\end{theo}
The proof of Theorem \ref{thm:eigenvectors_R} uses the invariances $x\to x+C$,  $\psi\to \psi +C$ of  \eqref{eq:main_with_t} and the fact that the Stokes wave family \eqref{eq:solution_with_P} depends on two parameters $a$ and $P$. Moreover, $R_3$ is an eigenvector if and only if $\p_a c|_{a=\eps}= 0$. 

For long-wave perturbations, we seek solutions to the linearized system (\ref{eq:v1})-(\ref{eq:v2}) of the form 
\begin{equation}
\label{eq:v12_vstar12}
    v_j(x,t)=e^{\lambda t}e^{i\mu x}v_j^*(x),
\end{equation}
where $v_j^*(x)$ are $2\pi$-periodic. We assume 
\bq\label{choice:mu} 
\mu=\frac{m_0}{n_0}\in (-\mez, \mez)\setminus\{0\}~\text{ is a small rational number and choose}~L=n_0 2\pi \gg 2\pi,
\eq
 so that $v_j(\cdot, t)$ is $L$-periodic. The parameter $\mu$ is known as the Floquet exponent.

%
\subsection{Hamiltonian formulation}
A linear system $\p_tf=Lf$ is called Hamiltonian if  $L=JB$, where
    \[
J:=\begin{bmatrix}
    0& \mathrm{Id}\\
    -\mathrm{Id}&0
\end{bmatrix}
\]
is the canonical symplectic matrix, and $B$ is self-adjoint. An operator  $T$ is called symplectic if $T^*JT=J$. Conjugation by symplectic transformations preserves the Hamiltonian structure.

With nonzero vorticity,  the water wave system \eqref{eq:main_with_t} is not Hamiltonian with respect to the canonical symplectic matrix $J$ \cite{ConstantinIvanovProdanov2008}. Nevertheless,  Wahl\'en  \cite{Wahlen_Hamiltonian} discovered that the linear change of variables  
\[
\xi=\psi-\frac{\gamma}{2}\pinv\Pi_0^\perp \eta,
\]
put  \eqref{eq:main_with_t} in the canonical Hamiltonian form for  $(\eta,\xi)$. Consequently, linearizing \eqref{eq:main_with_t} around the Stokes wave leads to a linear Hamiltonian system for $(\ol\eta,\ol\xi)$, where 
\begin{equation}
    \ol\xi=\ol\psi-\frac{\gamma}{2}\pinv\Pi_0^\perp \ol\eta.
\end{equation}
Since the  good-unknown transformation
\begin{equation}
\label{eq:l12_eta_xi}
    \vet{l_1}{l_2}=Z\vet{\ol\eta}{\ol\xi}:=
    \begin{bmatrix}
        1&0\\
        -B^* & 1
    \end{bmatrix}
    \vet{\ol\eta}{\ol\xi}
\end{equation}
is symplectic, $(l_1,l_2)$ satisfies a Hamiltonian system. We have the relations
\begin{align}
    \label{eq:v1_l1_temp}l_1&=v_1,\\
    \nonumber l_2&=\ol\xi-B^*\ol\eta=\ol\psi-\frac{\gamma}{2}\pinv\Pi_0^\perp \ol\eta-B^*\ol\eta\\
    \label{eq:v2_l2_temp}&=v_2-\frac{\gamma}{2}\pinv\Pi_0^\perp v_1.
\end{align}
To handle functions of the form \eqref{eq:v12_vstar12}, we define 
\[
S(\T_L)=\left\{f\in L^2_x(\T_L):~e^{-i\mu x}f(x)\text{~ is~} 2\pi-\text{periodic}\right\},
\]
 where we recall that $\mu=\frac{m_0}{n_0}\in (-\mez, \mez)\setminus\{0\}$ and $L=2\pi n_0$. In particular,  $v_j(\cdot , t)\in S(\T_L)$ for each $t$. We gather properties of $S(\T_L)$ in the following lemma. 
\begin{lemm}\label{lemm:Sclass}\cite[Lemma 2.6]{ChenNguyen}
\label{lemma: STL} Let $f(x)\in S(\T_L)$. Then 

(1) $\Pi_0 f=0$. 

(2) If  $a(D_L)f\in L^2(\T_L)$, then $a(D_L)f\in S(\T_L)$. In particular, $\mathcal{H}_L\tanh\big(\wt h|D_L|\big)f\in S(\T_L)$ for all $\wt h>0$, and  $\pinv f\in S(\T_L)$ if $\Pi_0 f=0$.
\end{lemm}
Since  $v_1\in S(\T_L)$, Lemma \ref{lemm:Sclass} (1) implies  $\Pi_0^\perp v_1=v_1$, so that the relation between $(l_1, l_2)$  and $(v_1, v_2)$ simplifies to
\begin{align}    \label{eq:vjlj}
l_1&=v_1,\quad l_2=v_2-\frac{\gamma}{2}\pinv v_1.
\end{align}
\subsection{Conformal mapping}
Following \cite{Huy21}, we use a conformal mapping to flatten the fluid domain $\{-h<y<\eta^*(x)\}$. There exist a ``conformal depth correction'' $\fe$ and  a holomorphic bijection map $z_1(x, y)+iz_2(x, y)$ from $\{-(h+\fe)<y<0\}$ onto  $\{-h<y<\eta^*(x)\}$. Since $\eta^*$ is even and $2\pi$-periodic, $z_1$ is odd in $x$, $z_2$ is even in $x$, and both $z_1(x, y)-x$ and $z_2(x, y)$ are $2\pi$-periodic in $x$. We define the Riemann stretch 
\begin{equation}
    \zeta(x)=z_1(x,0).
\end{equation}
 $\zeta$ and  $\fe$ are determined as a fixed point of the system (see \cite[(2.14)]{Berti_finite_depth})
\begin{equation}
\label{eq:fixed_point_zeta_fe}
    \zeta(x)=x+\frac{\mathcal{H}}{\tanh\big((h+\fe)|D|\big)}\Pi_0^\perp\zeta_\sharp(\eta^*),\quad \fe=\fint \zeta_\sharp(\eta^*)dx.
\end{equation}
It follows from the analyticity of $\eps\mapsto \eta^*$ and the analytic implicit function theorem that  $\zeta$ and $f_\eps$ are analytic in $\eps$.  Moreover, using \eqref{eq:fixed_point_zeta_fe} and the expansion for $\eta^*$ in \eqref{expand:Stokes_without_g_all_depth}, we obtain 
\begin{lemm}\label{lemm:zeta} There hold the expansions
\begin{equation} 
    \label{eq:expand_zeta}
    \zeta(x)=x+\eps\sin(x)+\frac{\eps^2}{4}\Big(\frac{1}{\th}+\frac{c_0^2+\th+c_0^2\th ^2}{c_0^4}\Big)\frac{1+\th^2}{2\th}\sin(2x)+\rem 3,
\end{equation}
\begin{equation}
    \label{eq:expand_fe}
    \fe=\frac{\eps^2}{4}(c_0^2-4\th +c_0^2\th^2)+\rem 3.
\end{equation}
\end{lemm}
\begin{proof}
See Appendix \ref{appendix:expansion_conformal_mapping}.
\end{proof}
Next, we define the operators 
\begin{equation}
    \zeta_\sharp f=f\circ \zeta ,  \qquad   \zeta_*f=\zeta' (f\circ \zeta),
\end{equation}
so that  $\zeta_*\p_xf=\p_x(\zeta_\sharp f)$ for $f:\R\to \R$. From \cite[(2.13)]{Berti_finite_depth},  the Dirichlet-Neumann operator  $G(\eta^*)$ can be written as 
\begin{equation}
    G\left(\eta^*\right) f=\partial_x\left(\zeta_{\sharp}^{-1} \mathcal{H}_L\tanh\big((h+\fe)|D_L|\big)\zeta_{\sharp} f\right),
\end{equation}
where $\mathcal H_Lu$ is the Hilbert transform
\begin{equation}
    \widehat{\mathcal{H}_L u}^L(k)=-i \operatorname{sign}(k) \widehat{u}^L(k) .
\end{equation}
Further properties of $\zeta_\sharp$ and $\zeta_*$ are given in the next lemma.
\begin{lemm}\cite[Lemma 2.7]{ChenNguyen} \label{lemma:zeta}
   (1) For $f$, $g\in L^2(\T_L)$, there hold the identities 
    \begin{equation}
        \Pi_0 \zeta_*f=\Pi_0 f
    \end{equation}
  and
    \begin{equation}
        \pinv\zeta_*g=\Pi_0^\perp \zeta_\sharp\pinv g\quad\text{if~}\Pi_0 g=0. 
    \end{equation}
   (2) If $f\in S(\T_L)$, then  $\zeta_*f,\zeta_\sharp f\in S(\T_L)$. Conversely, if $g$ satisfies either $\zeta_*g\in S(\T_L)$ or $\zeta_\sharp g\in S(\T_L)$, then $g\in S(\T_L)$. 
\end{lemm}
We make the change of variables 
\bq\label{eq:ljwj} 
 w_1=\zeta_* l_1,\quad  w_2=\zeta_\sharp l_2.
\eq
It can be verified directly that the preceding transformation is symplectic, so that $(w_1, w_2)$ again satisfies a Hamiltonian system. Combining the relations \eqref{eq:vjlj} and \eqref{eq:ljwj}, we find the following relation between $v_j$ and $w_j$:
\bq\label{eq:vjwj}
   w_1=\zeta_* v_1,\quad w_2=\zeta_\sharp v_2-\frac{\gamma}{2}\pinv \zeta_* v_1,
\eq
We have used  Lemmas \ref{lemma: STL} and \ref{lemma:zeta} to have
\begin{equation*}
    v_1\in S(\T_L)\implies \zeta_\sharp \pinv v_1\in S(\T_L)\implies \zeta_\sharp \pinv v_1=\Pi_0^\perp \zeta_\sharp \pinv v_1=\pinv \zeta_* v_1.
\end{equation*}
We set 
\begin{equation}
    p=\frac{c^*-\zeta_\sharp V^*+\gamma\zeta_\sharp \eta^*}{\zeta'},\qquad q=-p\p_x(\zeta_\sharp B^*).
\end{equation}
For any functions $f_1,f_2\in H^1(\T_L)$, noting $\zeta_\sharp\Pi_0=\Pi_0$, there hold the identities
\bq
\label{identity:zeta:1}
    \begin{aligned}
        & \zeta_*\Big(\p_x\big((c^*-V^*+\gamma\eta^*)f_1\big)+G(\eta^*)f_2\Big)=\p_x (p(x)\zeta_*f_1)+|D_L|\tanh\big((h+\fe)|D_L|\big)(\zeta_\sharp f_2)
        \end{aligned}
        \eq
        and
\bq\label{identity:zeta:2}
\begin{aligned}
        & \zeta_\sharp\Big(-\big(1+(V^*-c^*-\gamma\eta^*)\p_xB^*\big)f_1+\gamma\Pi_0(V^*f_1)+(c^*-V^*+\gamma\eta^*)\p_xf_2+\gamma\p_x^{-1}G(\eta^*)f_2\Big)\\
        &\quad=-\frac{1+q(x)}{\zeta'(x)}\zeta_*f_1+\gamma\Pi_0(V^*f_1)+p(x)\p_x(\zeta_\sharp f_2)\\
        &\qquad+\gamma\Big(\mathcal{H}_L\tanh\big((h+\fe)|D_L|\big)(\zeta_\sharp f_2)-\Pi_0\zeta_\sharp^{-1}\mathcal{H}_L\tanh\big((h+\fe)|D_L|\big)(\zeta_\sharp f_2)\Big).
    \end{aligned}
\eq
We apply $\zeta_*$ to \eqref{eq:v1} and $\zeta_\sharp$ to \eqref{eq:v2}, appeal to the identities  \eqref{identity:zeta:1}-\eqref{identity:zeta:2}, and use  Lemmas \ref{lemma: STL} and \ref{lemma:zeta}  to have
\[
V^*v_1\in S(\T_L)\implies \Pi_0(V^* v_1)=0
\]
 and 
 \[
 \zeta_\sharp^{-1}\mathcal{H}_L\tanh\big((h+\fe)|D_L|\big)(\zeta_\sharp v_2)\in S(\T_L)\implies \Pi_0\zeta_\sharp^{-1}\mathcal{H}_L\tanh\big((h+\fe)|D_L|\big)(\zeta_\sharp v_2)=0.
 \]
We obtain
\begin{equation*}
    \begin{aligned}
        \p_t \zeta_*v_1
        &=\p_x (p(x)\zeta_*v_1)+|D_L|\tanh\big((h+\fe)|D_L|\big)(\zeta_\sharp v_2),\\
        \p_t\zeta_\sharp v_2 
        &=-\frac{1+q(x)}{\zeta'(x)}\zeta_*v_1+p(x)\p_x(\zeta_\sharp v_2)+\gamma\mathcal{H}_L\tanh\big((h+\fe)|D_L|\big)(\zeta_\sharp v_2).
    \end{aligned}
\end{equation*}
Combining this and the relation \eqref{eq:vjwj}, we obtain 
\begin{align*}
    \p_t w_1&=\p_x(pw_1)+|D_L|\tanh\big((h+\fe)|D_L|\big)\Big(w_2+\frac{\gamma}{2}\pinv w_1\Big),
\end{align*}
\begin{align*}
    \p_t \Big(w_2+\frac{\gamma}{2}\pinv w_1\Big)&=-\frac{1+q}{\zeta'}w_1+p\p_x(\zeta_\sharp v_2)+\gamma\mathcal{H}_L\tanh\big((h+\fe)|D_L|\big)(\zeta_\sharp v_2).
\end{align*}
It follows that 
\begin{align*}
    \p_t w_2 &=-\frac{1+q}{\zeta'}w_1+p\p_x \Big(w_2+\frac{\gamma}{2}\pinv w_1\Big)+\gamma \mathcal{H}_L\tanh\big((h+\fe)|D_L|\big)\Big(w_2+\frac{\gamma}{2}\pinv w_1\Big)-\frac{\gamma}{2}\pinv\p_tw_1\\
    &=-\frac{1+q}{\zeta'}w_1+p\p_xw_2+p\frac{\gamma}{2}w_1+\gamma \mathcal{H}_L\tanh\big((h+\fe)|D_L|\big)\Big(w_2+\frac{\gamma}{2}\pinv w_1\Big)\\
    &\quad-\frac{\gamma}{2}\Bigg(pw_1+\mathcal{H}_L\tanh\big((h+\fe)|D_L|\big)\Big(w_2+\frac{\gamma}{2}\pinv w_1\Big)\Bigg)\\
    &=-\frac{1+q}{\zeta'}w_1+p\p_xw_2+\frac{\gamma}{2}\mathcal{H}_L\tanh\big((h+\fe)|D_L|\big)\Big(w_2+\frac{\gamma}{2}\pinv w_1\Big).
\end{align*}
We have used the fact that $pw_1\in S(\T_L)$, giving  $\pinv \p_x(pw_1)=\Pi_0^\perp(p w_1)=pw_1$. We conclude that $(w_1, w_2)$ satisfies the system
\begin{equation}\label{sys:wj}
    \begin{aligned}
        \p_t\begin{bmatrix}
            w_1\\
            w_2
        \end{bmatrix}=\begin{bmatrix}
            \p_x\circ p+\frac{\gamma}{2}\mathcal{H}_L\tanh\big((h+\fe)|D_L|\big) & |D_L|\tanh\big((h+\fe)|D_L|\big)\\
            -\frac{1+q(x)}{\zeta^{\prime}(x)}-\frac{\gamma^2}{4}\frac{\tanh\big((h+\fe)|D_L|\big)}{|D_L|}& p\p_x+\frac{\gamma}{2}\mathcal{H}_L\tanh\big((h+\fe)|D_L|\big)
        \end{bmatrix}\begin{bmatrix}
            w_1\\
            w_2
        \end{bmatrix},
    \end{aligned}
\end{equation}
which is Hamiltonian. In the companion Mathematica notebook \verb|1 Riemann_stretch.nb|,  we expand the variable coefficients in \eqref{sys:wj}. The result is  recorded below. 
\begin{lemm} There hold the expansions
\[
    B^*=\eps c_0   \th \sin (x)+\eps ^2 \sin (2 x) \left[\left(\frac{1}{2 c_0^3}-c_0\right) \th+\frac{c_0}{2 \th}+\frac{1}{2 c_0}\right]+\rem 3,
    \]
    \[
    V^*=\eps c_0   \cos (x)+\eps ^2 \left\{\cos (2 x) \left[\left(\frac{1}{4 c_0^3}-\frac{c_0}{4}\right) \th^2+\frac{1}{4 c_0^3}+\frac{\th}{4 c_0}+\frac{1}{4 c_0 \th}+\frac{c_0}{4 \th^2}-\frac{c_0}{2}\right]+\frac{c_0 \th^2}{2}\right\}+\rem 3,
\]
\begin{align*}
    p&=c_0-\frac{\eps  \cos (x) \left(c_0^2+\th\right)}{c_0}+\frac{\eps ^2 }{8 c_0^5 \th^2 \left(c_0^2+\th\right)}\Bigg[c_0^8+\left(2 c_0^{10}+2 c_0^6\right) \th+\left(7 c_0^8+5 c_0^4+1\right) \th^4\\
    &\quad +\left(3 c_0^4-5 c_0^8\right) \th^2+\left(8 c_0^6+2 c_0^2\right) \th^3+\left(-2 c_0^{10}-2 c_0^6+2 c_0^2\right) \th^5\\
    &\quad +\cos (2 x) \Big(-2 c_0^8-6 c_0^6 \th-2 c_0^4 \th^6-6 c_0^2 \th^5+\left(2 c_0^8-6 c_0^4-2\right) \th^4\\
    &\quad+\left(2 c_0^8-8 c_0^4\right) \th^2+\left(2 c_0^6-6 c_0^2\right) \th^3\Big)\Bigg]+\rem 3,
    \end{align*}
    and
    \begin{align*}
    q&=-\eps c_0^2   \th \cos (x)+\eps ^2 \left[\cos (2 x) \left(-\frac{c_0^2}{\th}+\left(\frac{3 c_0^2}{2}-\frac{1}{c_0^2}\right) \th+\frac{\th^2}{2}-1\right)+\frac{c_0^2 \th}{2}+\frac{\th^2}{2}\right]+\rem 3\\
    \frac{1+q}{\zeta'}&=1+\eps  \cos (x) \left(-\th c_0^2 -1\right)+\eps ^2 \Bigg\{\cos (2 x) \Bigg[\left(\frac{1}{2}-\frac{1}{4 c_0^4}\right) \th^2-\frac{1}{4 c_0^4}-\frac{\th^3}{4 c_0^2}\\
    &\quad+\left(2 c_0^2-\frac{3}{2 c_0^2}\right) \th+\frac{-c_0^2-\frac{1}{4 c_0^2}}{\th}-\frac{1}{4 \th^2}-\frac{3}{4}\Bigg]+\frac{\th^2}{2}+c_0^2 \th+\frac{1}{2}\Bigg\}+\rem 3.
\end{align*}
\end{lemm}
\subsection{The linearized operator}
From the above transformations and the time dependence \eqref{eq:v12_vstar12} of $v_j$, we can set
\begin{align}
w_j&=e^{\lambda t}e^{i\mu x}u_j(x),
    \end{align}
where $u_j(x)$ is $2\pi$-periodic. In view of \eqref{sys:wj}, we arrive at the eigenvalue problem 
\begin{equation}
    \lambda U=\L_{\mu,\eps}U,\quad U=\vet{u_1}{u_2},
\end{equation}
where 
\begin{equation}\label{def:finalL}
    \L_{\mu,\eps}:=\begin{bmatrix}
        (i\mu+\p_x)\circ p-i\frac{\gamma}{2}\tanh\big((h+\fe)(D+\mu)\big)
        & (D+\mu)\tanh\big((h+\fe)(D+\mu)\big)\\
        -\frac{1+q}{\zeta'}-\frac{\gamma^2}{4}\frac{\tanh\big((h+\fe)(D+\mu)\big)}{D+\mu}
        & p(i\mu+\p_x)-i\frac{\gamma}{2}\tanh\big((h+\fe)(D+\mu)\big)
    \end{bmatrix}
\end{equation}
    and $\mu=\frac{m_0}{n_0}\in (-\mez, \mez)\setminus\{0\}$ as in \eqref{choice:mu}. We notice, however, that the preceding formula for $ \L_{\mu,\eps}$ is well-defined for all real numbers $\mu\in (-\mez, \mez)$, giving the natural extension of $\L_{\mu, \eps}: (H^1(\T))^2\to (L^2(\T))^2$.  The (extended) operator $\L_{\mu,\eps}$ is  Hamiltonian:
\begin{align*}
    \L_{\mu,\eps}&=J\B_{\mu,\eps},\\
    \B_{\mu,\eps}&=
    \begin{bmatrix}
        \frac{1+q}{\zeta'}+\frac{\gamma^2}{4}\frac{\tanh\big((h+\fe)(D+\mu)\big)}{D+\mu}
        & -p(i\mu+\p_x)+i\frac{\gamma}{2}\tanh\big((h+\fe)(D+\mu)\big)\\
        (i\mu+\p_x)\circ p-i\frac{\gamma}{2}\tanh\big((h+\fe)(D+\mu)\big)
        & (D+\mu)\tanh\big((h+\fe)(D+\mu)\big)
    \end{bmatrix},
\end{align*}
where $\B_{\mu,\eps}$ is self-adjoint. 

Since $\tanh(z)/z$ is analytic in $\Rr$,  the matrix operator $\L_{\mu,\eps}$ is clearly jointly analytic in $\mu$ and $\eps$ for small $\eps$. The analyticity in $\mu$ is due to finite depth. Indeed, as shown in our previous work \cite{ChenNguyen} for the infinite depth case, the $(2,1)$-entry of $  \L_{\mu,\eps}$ is replaced by 
\[
   -\frac{1+q}{\zeta'}-\frac{\gamma^2}{4}\frac{\sign(D+\mu)}{D+\mu}
   \]
   which is not even bounded as $\mu\to 0$. 
 
Moreover, $ \L_{\mu,\eps}$  is reversible in the sense that
    \begin{equation}
    \label{eq:def_rev}
        \cR\L_{\mu,\eps}=-\L_{\mu,\eps}\cR,
    \end{equation}
    where $\cR$ is the reversal operator 
    \begin{equation}
    \label{eq:def_rev_R}
        \cR\begin{bmatrix}
            u_1(x)\\
            u_2(x)
        \end{bmatrix}:=\begin{bmatrix}
            -\ol u_1(-x)\\
            \ol u_2(-x)
        \end{bmatrix}.
    \end{equation}
    Our goal in the subsequent sections is to fully describe the spectrum of $\L_{\mu,\eps}$ in a neighborhood of $0\in \Cc$. This spectrum depends on the constant vorticity $\gamma$ via the coefficients of $\L_{\mu,\eps}$.
    \section{Kato's perturbed basis}
\subsection{Eigenvalues of $\L_{0,\eps}$}\label{subsec:upev}
In this subsection, we discuss eigenvalues of $\L_{\mu,\eps}$ near $0\in \Cc$ when $\mu=0$.  When both $\mu$ and $\eps$ are zero, 
\begin{equation}
    \label{eq:L00_op}
    \L_{0,0}=
\begin{bmatrix}
    c_0\p_x-i\frac{\gamma}{2}\tanh(hD) & D\tanh(hD)\\
    -1-\frac{\gamma^2}{4}\frac{\tanh(hD)}{D} & c_0\p_x-i\frac{\gamma}{2}\tanh(hD)
\end{bmatrix}.
\end{equation}
The characteristic equation is
\begin{equation}
\label{eq:lambda}
\begin{aligned}
    0&=\det\Bigg(
    \begin{bmatrix}
    \lambda-c_0ik+i\frac{\gamma}{2}\tanh(hk) & -k\tanh(hk)\\
    1+\frac{\gamma^2}{4}\frac{\tanh(hk)}{k} & \lambda-c_0ik+i\frac{\gamma}{2}\tanh(hk)
\end{bmatrix}
    \Bigg)\\
    &=\Big(\lambda-c_0ik+i\frac{\gamma}{2}\tanh(hk)\Big)^2+\tanh(hk)\Big(k+\frac{\gamma^2}{4}\tanh(hk)\Big)\\
    &=\lambda^2+i\lambda\Bigg(-2c_0 k+\tanh(hk)\Big(\frac{c_0}{\th}-\frac{1}{c_0}\Big)\Bigg)-c_0^2\Big(k^2-k\frac{\tanh(hk)}{\th}\Big).
\end{aligned}
\end{equation}
Hence, the eigenvalues of $ \L_{0,0}$ are $\lambda=i\omega_{k}^\pm$, 
\begin{equation}
\label{eq:omega_k_pm}
    \omega_k^\pm = c_0k-\frac{\gamma}{2}\tanh(hk)\pm\sqrt{\tanh(hk)\Big(k+\frac{\gamma^2}{4}\tanh(hk)\Big)}.
\end{equation}
Since 
\[
\sqrt{\tanh(h)\Big(1+\frac{\gamma^2}{4}\tanh(h)\Big)}=\sqrt{\th+\frac{1}{4}\th^2\Big(\frac{c_0}{\th}-\frac{1}{c_0}\Big)^2}=\mez \Big(c_0+\frac{\th}{c_0}\Big),
\]
we have
\begin{align*}
    \om_0^+&=0,&\om_0^-&=0,\\
    \om_1^+&=c_0+\frac{\th}{c_0},&\om_1^-&=0,\\
    \om_{-1}^+&=0,&\om_{-1}^-&=-(c_0+\frac{\th}{c_0}).
\end{align*}
Since $\tanh(x)$ is strictly concave on $(0, \infty)$, for $k\ge 2$, we have 
\[
\tanh(h)=\tanh\big(\frac{1}{k}hk+(1-\frac{1}{k})0\big)>\frac{1}{k}\tanh(hk).
\]
Hence, $k-\frac{\tanh(hk)}{\th}\ne 0$ for $|k|\ge 2$, proving that $\lambda=0$  is not a root of \eqref{eq:lambda}. In addition, $\om_k^\pm\sim c_0k$ when $|k|\to\infty$. We conclude that the spectrum of $\L_{0,0}$ is separated into two isolated parts $\sigma'(\L_{0,0})\cup\sigma''(\L_{0,0})$, where 
\bq\label{spectrum:00}
    \sigma'(\L_{0,0})=\{i\omega_{0}^+, i\omega_{0}^-, i\omega_{1}^-, i\omega_{-1}^+\},
\eq
all of which are $0$. By semicontinuity of  separated parts of the spectrum \cite[Theorem 3.16]{Kato}, for small $\mu$ and $\eps$, the spectrum of $\L_{\mu, \eps}$ is again separated into two parts, where the part near $0\in \Cc$ has a total dimension of generalized eigenspace of $4$. 

We recall that the operator $\L_{\mu, \eps}$ in \eqref{def:finalL} was obtained from $\Ls_\eps$ \eqref{eq:v_matrix} by the transformations \eqref{eq:vjlj} and \eqref{eq:ljwj}, which are defined for $\mu\ne 0$.  In the next lemma, we directly verify the relation between $\Ls_\eps$ and $\L_{0, \eps}$, where the latter was defined by extension of \eqref{def:finalL}. This relation will be used in Proposition \ref{prop:evL0eps} to obtain eigenvectors of $\L_{0, \eps}$ from those of $\Ls_\eps$.
\begin{lemm}   \label{prop:convert_RU} 
We define 
    \begin{equation}
        \Pi_{01}\begin{bmatrix}
            r_1\\
            r_2
        \end{bmatrix}=\Pi_0r_1
    \end{equation}
    and 
    \begin{equation}
\label{def:B0}
    B_0:=
    \begin{bmatrix}
        \zeta_* & 0\\
        -\frac{\gamma}{2}\zeta_\sharp\pinv & \zeta_\sharp
    \end{bmatrix}.
\end{equation}
For any $R\in (H^1(\T))^2$ satisfying $\Pi_{01} R= 0$, $B_0 R$ is well-defined and
\bq\label{relation2L}
B_0 \Ls_\eps  R= \L_{0,\eps} B_0 R+\vet{0}{s},
\eq
where  $s\in\C$ is a constant depending on $R$.
\end{lemm}
\begin{proof}
We denote $R=\vet{r_1}{r_2}$, where $\Pi_0r_1=0$. Then $\p_x^{-1}r_1$ is well-defined, and hence so is $B_0R$. We proceed to compute  $B_0\Ls_\eps  R$, where 
\[
\Ls_\eps R=
\vet{\p_x\big((c^*-V^*+\gamma\eta^*)r_1\big)+G(\eta^*)r_2}
{-\big(1+(V^*-c^*-\gamma\eta^*)\p_xB^*\big)r_1+\gamma\Pi_0(V^*r_1)+(c^*-V^*+\gamma\eta^*)\p_xr_2+\gamma\p_x^{-1}G(\eta^*)r_2}.
\]
By Lemma \ref{lemma:zeta} (1), if $\Pi_0f=0$, then 
    \[
    \zeta_\sharp\pinv  f= \Pi_0\zeta_\sharp\pinv f+\Pi_0^\perp \zeta_\sharp\pinv f=\Pi_0\zeta_\sharp\pinv f+\pinv\zeta_*f,
    \]
implying that 
    \[
     B_0=
    \begin{bmatrix}
        \zeta_* & 0\\
        -\frac{\gamma}{2}\pinv\zeta_*  -\frac{\gamma}{2}\Pi_0\zeta_\sharp\pinv& \zeta_\sharp
    \end{bmatrix}
    \]
when acting on vectors whose first component has mean zero. Using this and   \eqref{identity:zeta:1}-\eqref{identity:zeta:2},  we compute
\begin{align*}
(B_0\Ls_\eps  R)_1&=\p_x(p\zeta_*r_1)+|D|\tanh\big((h+\fe)|D|\big)(\zeta_\sharp r_2)
\end{align*}
and 
\begin{align*}
(B_0\Ls_\eps  R)_2 &=
-\frac{\gamma}{2}\Big(\Pi_0^\perp(p\zeta_*r_1)+\mathcal{H}\tanh\big((h+\fe)|D|\big)(\zeta_\sharp r_2)\Big)\\
&\qquad -\frac{\gamma}{2}\Pi_0\zeta_\sharp\Big(\Pi_0^\perp\big((c^*-V^*+\gamma\eta^*)r_1\big)+\pinv G(\eta^*)r_2\Big)\\
      & \qquad -\frac{1+q}{\zeta'}\zeta_*r_1+\gamma\Pi_0(V^*r_1)+p\p_x(\zeta_\sharp r_2)\\
       & \qquad +\gamma\Big(\mathcal{H}\tanh\big((h+\fe)|D|\big)(\zeta_\sharp r_2)-\Pi_0\zeta_\sharp^{-1}\mathcal{H}\tanh\big((h+\fe)|D|\big)(\zeta_\sharp r_2)\Big)\\
       &=-\frac{\gamma}{2}p\zeta_*r_1 -\frac{1+q}{\zeta'}\zeta_*r_1+p\p_x(\zeta_\sharp r_2)+\frac{\gamma}{2} \mathcal{H}\tanh\big((h+\fe)|D|\big)(\zeta_\sharp r_2)+s,
        \end{align*}
        where
        \begin{align*}
        s&:=\frac{\gamma}{2}\Pi_0(p\zeta_*r_1)+\gamma\Pi_0(V^*r_1)\\
        &\qquad -\frac{\gamma}{2}\Pi_0\zeta_\sharp\Big(\Pi_0^\perp\big((c^*-V^*+\gamma\eta^*)r_1\big)+\pinv G(\eta^*)r_2\Big) -\gamma\Pi_0\zeta_\sharp^{-1}\mathcal{H}\tanh\big((h+\fe)|D|\big)(\zeta_\sharp r_2)\in \Cc.
        \end{align*}
We recall that $\L_{0, \eps}$ was defined by extending \eqref{def:finalL} to $\mu=0$, that is,
\bq\label{form:L0eps}
\L_{0, \eps}=  \begin{bmatrix}
        \p_x\circ p-i\frac{\gamma}{2}\tanh\big((h+\fe)D\big)
        & D\tanh\big((h+\fe)D\big)\\
        -\frac{1+q}{\zeta'}-\frac{\gamma^2}{4}\frac{\tanh\big((h+\fe)D\big)}{D}
        & p\p_x-i\frac{\gamma}{2}\tanh\big((h+\fe)D\big)
    \end{bmatrix}
\eq
Hence, 
\begin{align*}
(\L_{0,\eps} B_0R)_1&= \p_x(p\zeta_* r_1)-i\frac{\gamma}{2}\tanh\big((h+\fe)D\big)\zeta_* r_1-\frac{\gamma}{2}D\tanh\big((h+\fe)D\big)\pinv \zeta_* r_1\\
&\qquad+D\tanh\big((h+\fe)D\big)\zeta_\sharp r_2\\
&=\p_x(p\zeta_* r_1)+D\tanh\big((h+\fe)D\big)\zeta_\sharp r_2
\end{align*}
since 
\[
D\tanh\big((h+\fe)D\big)\pinv \zeta_* r_1=-i\tanh\big((h+\fe)D\big)\p_x\pinv \zeta_* r_1=-i\tanh\big((h+\fe)D\big)\zeta_* r_1.
\]
On the other hand,
\begin{align*}
(\L_{0,\eps} B_0R)_2&=-\frac{1+q}{\zeta'}\zeta_* r_1-\frac{\gamma^2}{4}\frac{\tanh\big((h+\fe)D\big)}{D}\zeta_* r_1\\
&\qquad-\frac{\gamma}{2}p\Pi_0^\perp \zeta_* r_1+i\frac{\gamma^2}{4}\tanh\big((h+\fe)D\big)\pinv\zeta_* r_1+p\p_x\zeta_\sharp r_2-i\frac{\gamma}{2}\tanh\big((h+\fe)D\big)\zeta_\sharp r_2\\
&=-\frac{1+q}{\zeta'}\zeta_* r_1-\frac{\gamma}{2}p \zeta_* r_1+p\p_x\zeta_\sharp r_2-i\frac{\gamma}{2}\tanh\big((h+\fe)D\big)\zeta_\sharp r_2
\end{align*}
since $\Pi_0\zeta_* r_1=\Pi_0 r_1=0$ by Lemma \ref{lemma:zeta} (1).    

\eqref{relation2L} follows from the above computations. 
\end{proof}
\subsection{Kato's similarity transformation}
In order to track the spectrum of $\L_{\mu, \eps}$ as $(\mu, \eps)$ is perturbed from $(0, 0)$, we follow  \cite{Berti_Full} and use Kato's similarity transformation
\begin{equation}
\label{def:Kato}
    \cK_{\mu,\eps} := \{1-(P_{\mu,\eps}-P_{0, 0})^2\}^{-1/2} \{P_{\mu,\eps} P_{0, 0} + (1-P_{\mu,\eps})(1-P_{0, 0})\}.
\end{equation}
Here,  $P_{\mu,\eps}:  (L^2(\T))^2\to (H^1(\T))^2$ is the spectral projection
\begin{equation}
    \label{def:P_spec_proj}
P_{\mu,\eps} = -\frac1{2\pi i} \oint_\Gamma (\L_{\mu,\eps} - \lambda)^{-1} d\lambda,
\end{equation}
where $\Gamma$ is a closed circle of sufficiently small radius around $0$ such that $\Gamma\subset \rho(\L_{\mu,\eps})$ for $(\mu,\eps)$ close to $(0,0)$. The analyticity of $\L_{\mu,\eps}$ implies that the operators $\cK_{\mu,\eps}$ and $P_{\mu,\eps}$ are analytic in $(\mu,\eps)$ near $(0,0)$. 

We set  
\begin{equation}
    \cV_{\mu,\eps}:=\operatorname{Ran}(P_{\mu,\eps}),
\end{equation}
From \cite[Lemma 3.1]{Berti_Full}, we have the following properties 
\begin{align}
    \L_{\mu,\eps}:\cV_{\mu,\eps}&\to\cV_{\mu,\eps},\label{prop1_berti_kato}\\
    \sigma(\L_{\mu,\eps})\cap\{z\in\mathbb{C} \text{ inside }\Gamma\}&=\sigma(\L_{\mu,\eps}|_{\cV_{\mu,\eps}}),\label{prop2_berti_kato}\\
    \cV_{\mu,\eps}&=\cK_{\mu,\eps}\cV_{0,0}\label{prop3_berti_kato}
\end{align}
Therefore, the spectrum of $\L_{\mu, \eps}$ near $0\in \Cc$ can be determined from the action of $\L_{\mu, \eps}$ on the finite-dimensional space $\cV_{\mu,\eps}=\cK_{\mu,\eps}\cV_{0,0}$. To this end, we first fix a basis for $\cV_{0,0}$. At frequencies  $k=\pm 1$,  $\cL_{0, 0}$ is represented  by the matrix 
\[
\begin{bmatrix}
    \pm i\Big(c_0-\frac{\gamma}{2}\th\Big) & \th\\
    -1-\frac{\gamma^2}{4}\th &  \pm i\Big(c_0-\frac{\gamma}{2}\th\Big)
\end{bmatrix}.
\]
Two independent eigenvectors corresponding to the eigenvalue $0$ are $\vet{\th e^{\pm ix}}{\mp i(c_0-\frac{\gamma}{2}\th)e^{\pm ix}}$. Taking linear combinations, we get the linear eigenvectors
\begin{equation}
\label{eq:U1_pm}
U_{1,+}=\vet{\th \cos x}{\Big(c_0-\frac{\gamma}{2}\th\Big)\sin x},\quad U_{1,-}=\vet{-\th \sin x}{\Big(c_0-\frac{\gamma}{2}\th\Big)\cos x}.
\end{equation}
At frequency $k=0$, $\L_{0, 0}$ is given 
\[
\begin{bmatrix}
0 & 0\\
-1-\frac{\gamma^2}{4}h& 0
\end{bmatrix}.
\]
For the vectors 
\begin{equation}
\label{eq:U0_pm}
U_{0,+}=\vet{1}{0},\quad U_{0,-}=\vet{0}{1},
\end{equation}
we find  $\L_{0,0} U_{0,-}=0$ and $\L_{0,0} U_{0,+}=(-1-\frac{\gamma^2}{4}h)U_{0,-}$. Therefore, 
\bq
\cV_{0,0}=\Span\{U_{1,+},U_{1,-},U_{0,+},U_{0,-}\},\quad \cV_{\mu,\eps}=\Span \{U_{1,+}^{\mu,\eps},U_{1,-}^{\mu,\eps},U_{0,+}^{\mu,\eps},U_{0,-}^{\mu,\eps}\},
\end{equation}
where
\begin{equation}
\label{def:U_Katobasis}
U_{j,\pm}^{\mu,\eps}=\cK_{\mu,\eps} U_{j,\pm}.
\end{equation}
The vectors  $U_{j,\pm}^{\mu,\eps}$ are analytic in $(\mu,\eps)$ near $(0,0)$ because $\cK_{\mu,\eps}$ is so. 

Since $\L_{\mu, \eps}$ is a reversible Hamiltonian operator,  \cite[Lemma 3.2 (ii)]{Berti_Full} implies that  $\cK_{\mu,\eps}$ is symplectic and reversibility-preserving, i.e. $\cK_{\mu,\eps}^*J\cK_{\mu,\eps}=J$ and $\cK_{\mu, \eps}\cR=\cR\cK_{\mu, \eps}$. It follows that
\begin{equation}
\label{eq:symplectic_cK}
(J\cK_{\mu,\eps}U,\cK_{\mu,\eps}V)=(\cK_{\mu,\eps}^*J\cK_{\mu,\eps}U,V)=(JU,V),
\end{equation}
and hence 
\bq\label{eq:JUmn_general:0}
 (JU_{m,\sigma}^{\mu,\eps},U_{n,\tau}^{\mu,\eps})=
    (JU_{m,\sigma},U_{n,\tau}).
    \eq
Direct calculation gives
\begin{equation}
\begin{aligned}
&(JU_{1,+},U_{1,-})=-2\pi\th \Big(c_0-\frac{\gamma}{2}\th\Big)=-(JU_{1,-},U_{1,+}),\\
&(JU_{0,+},U_{0,-})=-2\pi=-(JU_{0,-},U_{0,+}),
\end{aligned}
\end{equation}
and all other pairs of the form  $(JU_{j,\sigma},U_{k,\tau})$ are zero. Consequently, denoting 
\begin{equation}
    \label{def:beta}
    \beta_1=\th\Big(c_0-\frac{\gamma}{2}\th\Big),\quad \beta_0=1,
\end{equation}
we deduce from \eqref{eq:JUmn_general:0} that
\begin{equation}
\label{eq:JUmn_general}
    (JU_{m,\sigma}^{\mu,\eps},U_{n,\tau}^{\mu,\eps})=
    (JU_{m,\sigma},U_{n,\tau})=-
    2\pi \delta_{mn}\delta_{\sigma,-\tau}\sigma\beta_m.
\end{equation}

\begin{prop}
\label{prop:structure_Lmat}
    (1) The $4\times 4$ matrix representing the linear operator $\L_{\mu,\eps}=J\B_{\mu,\eps}:\cV_{\mu,\eps}\to \cV_{\mu,\eps}$ with respect to the basis $\{U_{1,+}^{\mu,\eps},U_{1,-}^{\mu,\eps},U_{0,+}^{\mu,\eps},U_{0,-}^{\mu,\eps}\}$ (in this order) is
    \begin{equation}
    \label{def:Lmat}
    \Lmat_{\mu,\eps}
    :=
    \frac{1}{2\pi}
    \begin{bmatrix}
        \displaystyle
        \frac{
        \left(
            \mathcal{B}_{\mu,\eps}U_{1,+}^{\mu,\eps},
            U_{1,-}^{\mu,\eps}
        \right)
        }{\beta_1}
        &
        \displaystyle
        -\frac{
        \left(
            \mathcal{B}_{\mu,\eps}U_{1,+}^{\mu,\eps},
            U_{1,+}^{\mu,\eps}
        \right)
        }{\beta_1}
        &
        \left(
            \mathcal{B}_{\mu,\eps}U_{1,+}^{\mu,\eps},
            U_{0,-}^{\mu,\eps}
        \right)
        &
        -\left(
            \mathcal{B}_{\mu,\eps}U_{1,+}^{\mu,\eps},
            U_{0,+}^{\mu,\eps}
        \right)
        \\[3mm]
        \displaystyle
        \frac{
        \left(
            \mathcal{B}_{\mu,\eps}U_{1,-}^{\mu,\eps},
            U_{1,-}^{\mu,\eps}
        \right)
        }{\beta_1}
        &
        \displaystyle
        -\frac{
        \left(
            \mathcal{B}_{\mu,\eps}U_{1,-}^{\mu,\eps},
            U_{1,+}^{\mu,\eps}
        \right)
        }{\beta_1}
        &
        \left(
            \mathcal{B}_{\mu,\eps}U_{1,-}^{\mu,\eps},
            U_{0,-}^{\mu,\eps}
        \right)
        &
        -\left(
            \mathcal{B}_{\mu,\eps}U_{1,-}^{\mu,\eps},
            U_{0,+}^{\mu,\eps}
        \right)
        \\[3mm]
        \displaystyle
        \frac{
        \left(
            \mathcal{B}_{\mu,\eps}U_{0,+}^{\mu,\eps},
            U_{1,-}^{\mu,\eps}
        \right)
        }{\beta_1}
        &
        \displaystyle
        -\frac{
        \left(
            \mathcal{B}_{\mu,\eps}U_{0,+}^{\mu,\eps},
            U_{1,+}^{\mu,\eps}
        \right)
        }{\beta_1}
        &
        \left(
            \mathcal{B}_{\mu,\eps}U_{0,+}^{\mu,\eps},
            U_{0,-}^{\mu,\eps}
        \right)
        &
        -\left(
            \mathcal{B}_{\mu,\eps}U_{0,+}^{\mu,\eps},
            U_{0,+}^{\mu,\eps}
        \right)
        \\[3mm]
        \displaystyle
        \frac{
        \left(
            \mathcal{B}_{\mu,\eps}U_{0,-}^{\mu,\eps},
            U_{1,-}^{\mu,\eps}
        \right)
        }{\beta_1}
        &
        \displaystyle
        -\frac{
        \left(
            \mathcal{B}_{\mu,\eps}U_{0,-}^{\mu,\eps},
            U_{1,+}^{\mu,\eps}
        \right)
        }{\beta_1}
        &
        \left(
            \mathcal{B}_{\mu,\eps}U_{0,-}^{\mu,\eps},
            U_{0,-}^{\mu,\eps}
        \right)
        &
        -\left(
            \mathcal{B}_{\mu,\eps}U_{0,-}^{\mu,\eps},
            U_{0,+}^{\mu,\eps}
        \right)
    \end{bmatrix}.
\end{equation}
(2) The $(j, k)$-entry of $\Lmat_{\mu,\eps}$ is real (respectively, pure imaginary) if $j$ and $k$ have different (respectively, same)  parity.
\end{prop}
\begin{proof}
    (1)   By \ref{prop1_berti_kato}, the space $\cV_{\mu,\eps}$ is actually an invariant subspace of $\L_{\mu,\eps}$.

For any $U=\sum_{0\le j\le 1,\sigma=\pm}a_{j,\sigma}U_{j,\sigma}^{\mu,\eps}\in \cV_{\mu,\eps}$, \eqref{eq:JUmn_general} implies 
       \[
       (JU,U_{k,\tau}^{\mu,\eps})=a_{k,-\tau} 2\pi \tau \beta_k,\quad k=0,1, \quad \tau=\pm.
       \]
Consequently, 
       \begin{equation}
    \label{eq:expr_cV_vector_in_basis_by_inner_product}
        U= \sum_{0\le k\le 1,\tau=\pm} -\frac{1}{2\pi}\frac{\tau}{\beta_k}(JU,U_{k,-\tau}^{\mu,\eps})U_{k,\tau}^{\mu,\eps}\quad\forall U\in  \cV_{\mu,\eps}.
    \end{equation}
We recall from  \eqref{prop1_berti_kato} that $\cV_{\mu, \eps}$ is an invariant subspace of $\L_{\mu, \eps}$. Hence, applying \eqref{eq:expr_cV_vector_in_basis_by_inner_product} to  $U=\L_{\mu,\eps}U_{j,\sigma}^{\mu,\eps}$ gives
\begin{align*}
    \L_{\mu,\eps}U_{j,\sigma}^{\mu,\eps}
     &=\sum_{0\le k\le 1,\tau=\pm} -\frac{1}{2\pi}\frac{\tau}{\beta_k}(J\L_{\mu,\eps}U_{j,\sigma}^{\mu,\eps},U_{k,-\tau}^{\mu,\eps})U_{k,\tau}^{\mu,\eps}\\
     &=\sum_{0\le k\le 1,\tau=\pm} \frac{1}{2\pi}\frac{\tau}{\beta_k}(\B_{\mu,\eps} U_{j,\sigma}^{\mu,\eps},U_{k,-\tau}^{\mu,\eps})U_{k,\tau}^{\mu,\eps}.
\end{align*}
Therefore, the matrix representation of $\cL_{\mu, \eps}$ in $\cV_{\mu, \eps}$ is given by \eqref{def:Lmat}. 

      (2) From  the reversibility \eqref{eq:def_rev} of $\cL_{\mu, \eps}$ and the identity  $\cR J=-J\cR$, we have $\cR\B_{\mu,\eps}=\B_{\mu,\eps}\cR$. It is readily verified from the formulas \eqref{eq:U1_pm} and \eqref{eq:U0_pm} for $U_{j, \sigma}$ that  $\cR U_{j,\sigma}=-\sigma U_{j,\sigma}$.  Combining this with the fact that  $\cK_{\mu,\eps}$ is reversibility-preserving, we deduce 
\begin{equation}
    \label{eq:Uj_reversible}
        \cR U_{j,\sigma}^{\mu,\eps}=-\sigma U_{j,\sigma}^{\mu,\eps},
    \end{equation}
Invoking the identity $(\cR v,\cR w)=\ol{(v,w)}$, we find
 \[
 (\B_{\mu,\eps} U_{j,\sigma}^{\mu,\eps},U_{k,-\tau}^{\mu,\eps})=
 \ol{(\cR\B_{\mu,\eps} U_{j,\sigma}^{\mu,\eps}, \cR U_{k,-\tau}^{\mu,\eps})}=
 -\sigma\tau \ol{(\B_{\mu,\eps} U_{j,\sigma}^{\mu,\eps},U_{k,-\tau}^{\mu,\eps})}.
 \]
 Hence, $(\B_{\mu,\eps} U_{j,\sigma}^{\mu,\eps},U_{k,-\tau}^{\mu,\eps})$ is real (respectively, purely imaginary) if $\sigma \ne  \tau$ (respectively, $\sigma =\tau$).  Equivalently, the $(j, k)$-entry of $\Lmat_{\mu,\eps}$ is real (respectively, pure imaginary) if $j$ and $k$ have different (respectively, same)  parity.
 \end{proof}
\section{Analysis of the matrix $\Lmat_{\mu,\eps}$}
\subsection{Second-order expansion of $\Lmat_{\mu,\eps}$}
\def\errorTol{O(\eps^3,\mu\eps^2,\mu^2\eps,\mu^3)}
In this subsection, we  expand $\Lmat_{\mu,\eps}$  to second order in $\mu$ and $\eps$, leaving error terms $\errorTol$. In other words, we aim to find the coefficients $ \Lmat^{j,k}$ in the expansion
\begin{equation}
    \label{expansion:Lmat_formal}
    \Lmat_{\mu,\eps}=\sum_{j+k\le 2} \Lmat^{j,k}\mu^j\eps^k+\errorTol.
\end{equation}
Denote 
\begin{equation}
    \label{def:Delta} \Delta\L_{\mu,\eps}=\L_{\mu,\eps}-\L_{0,0},
    \quad \Delta P_{\mu,\eps}=P_{\mu,\eps}-P_{0,0}.
\end{equation}
To calculate the inverse operator in $P_{\mu, \eps}$ \eqref{def:P_spec_proj}, we use the Neumann series 
\begin{equation}
    \label{eq:Neumann_series}
    (\L_{\mu,\eps}-\lambda)^{-1}=(\L_{0,0}-\lambda)^{-1}
  +\sum_{k=1}^\infty (-1)^k(\L_{0,0}-\lambda)^{-1}\Big(\Delta \L_{\mu,\eps} (\L_{0,0}-\ld)^{-1}\Big)^k,
\end{equation}
where $\ld\in\Gamma\subset \rho(\cL_{\mu, \eps})$. Since $\L_{\mu,\eps}$ is analytic in $\mu$ and $\eps$, we have $\Delta\L_{\mu,\eps}=O(\mu,\eps)$ as an operator in $\cL((H^1)^2,(L^2)^2)$. Consequently,
\[
(-1)^k(\L_{0,0}-\lambda)^{-1}\Big(\Delta \L_{\mu,\eps} (\L_{0,0}-\ld)^{-1}\Big)^k=\errorTol\quad \text{for } k\ge 3.
\]
This implies
\begin{multline}
    \label{cutoff:Neumann}
    (\L_{\mu,\eps}-\lambda)^{-1}-(\L_{0,0}-\lambda)^{-1}=-(\L_{0,0}-\lambda)^{-1}\Delta\L_{\mu,\eps}(\L_{0,0}-\lambda)^{-1}\\
    +(\L_{0,0}-\lambda)^{-1}\Delta\L_{\mu,\eps}(\L_{0,0}-\lambda)^{-1}\Delta\L_{\mu,\eps}(\L_{0,0}-\lambda)^{-1}+\errorTol,
\end{multline}
and hence 
\begin{equation}
    \label{cutoff:P}
    \begin{aligned}
        (\Delta P_{\mu,\eps})U &=
        -\res \Big((\L_{\mu,\eps}-\lambda)^{-1}U-(\L_{0,0}-\lambda)^{-1} U\Big)\\
        &=-[\lambda^{-1}] \Big(-(\L_{0,0}-\lambda)^{-1}\Delta\L_{\mu,\eps}(\L_{0,0}-\lambda)^{-1}U\\
    &\quad +(\L_{0,0}-\lambda)^{-1}\Delta\L_{\mu,\eps}(\L_{0,0}-\lambda)^{-1}\Delta\L_{\mu,\eps}(\L_{0,0}-\lambda)^{-1}U\Big)+\errorTol.
    \end{aligned}
\end{equation}
The notation $[X_j^m]F(X_1, \dots, X_n)$ is defined by \eqref{nota:coefficient}. We recall from \eqref{eq:L00_op} that $\L_{0, 0}$ is the Fourier multiplier 
\begin{equation}
    \label{eq:L00_op_abc}
    \L_{0,0}=
    \begin{bmatrix}
        a(D) & c(D)\\
        b(D) & a(D)
    \end{bmatrix},
\end{equation}
\begin{equation}
\label{def:abc}
    a(k):=ic_0k-i\frac{\gamma}{2}\tanh(hk),\quad
    b(k):=-1-\frac{\gamma^2}{4}\frac{\tanh(hk)}{k},\quad
    c(k):=k\tanh(hk).
\end{equation}
Hence, we have
\begin{equation}
    \label{eq:Linv_op_abc}
    (\L_{0,0}-\lambda)^{-1}
    =\frac{1}{(a(D)-\lambda)^2-b(D)c(D)}\begin{bmatrix}
        a(D)-\lambda &-c(D)\\
       - b(D)& a(D)-\lambda
    \end{bmatrix}.
\end{equation}
\def\bracketLd#1{_{[#1]}}
We recall from \eqref{spectrum:00} that the generalized eigenspace for the eigenvalue $\ld=0$ of $\cL_{0, 0}$  consists of vectors having frequencies $k=0, \pm 1$. Hence, $a(k)^2-b(k) c(k)\ne 0$ for $|k|\ge 2$. For $p\ge 1$, we denote by $\Pi_{\ge p}$ the projection from $L^2(\T)$ onto $\Span\{e^{ikx}: |k|\ge p\}$. For $m\ge 0$ and $p\ge 1$, we set 
\[
H^m_{\ge p}(\T)=\left\{f\in H^m(\T):~\hat{f}(k)=0\quad\forall |k|<p \right\}, \quad m\ge 0.
\]
 Then,  $\L_{0,0}: (H^1)^2_{\ge 2}\to (L^2)^2_{\ge 2}$ is a homeomorphism. Hence, we have the Neumann series expansion 
\begin{equation}\label{invL00:1}
    (\L_{0,0}-\lambda)^{-1}\Pi_{\ge 2}=(\L_{0,0})^{-1}\Pi_{\ge 2}+\lambda(\L_{0,0})^{-2}\Pi_{\ge 2}+\lambda^2 (\L_{0,0})^{-3}\Pi_{\ge 2}+\cdots.
\end{equation}
At frequencies $k=\pm 1$, $a(\pm 1)=\pm(i c_0-i\frac{\gamma}{2}\th)=\pm\mez(ic_0+i\frac{\th}{c_0})\ne 0$, so the symbol of $ (\L_{0,0}-\lambda)^{-1}$ is 
\bq\label{invL00:2}
\begin{aligned}
   & \frac{1}{\lambda^2-2a(k)\lambda} \begin{bmatrix}
        a(k)-\lambda &-c(k)\\
       - b(k)& a(k)-\lambda
    \end{bmatrix}\\
    &=\sum_{j=-1}^{\infty} -(2a(k))^{-j-2}\lambda^j
     \begin{bmatrix}
        a(k)-\lambda &-c(k)\\
       - b(k)& a(k)-\lambda
    \end{bmatrix}\\
    &=\sum_{j=-1}^\infty  -(2a(k))^{-j-2}\lambda^j
     \begin{bmatrix}
        a(k) &-c(k)\\
       - b(k)& a(k)
    \end{bmatrix}
    +\sum_{j=0}^\infty (2a(k))^{-j-1}\lambda^j\mathrm{Id}_2.
\end{aligned}
\eq
For $k=0$, the symbol of $ (\L_{0,0}-\lambda)^{-1}$ is simply
\bq\label{invL00:3}
    \frac{1}{\lambda^2}
    \begin{bmatrix}
        -\lambda & 0\\
        1+\frac{\gamma^2}{4}h & -\lambda
    \end{bmatrix}.
\eq
From \eqref{invL00:1}, \eqref{invL00:2}, and \eqref{invL00:3}, we obtain
\begin{equation}\label{invL00:4}
    (\L_{0,0}-\lambda)^{-1}=\sum_{j=-2}^{\infty}\lambda^{j} R\bracketLd{j},
\end{equation}
Substituting \eqref{invL00:4} into \eqref{cutoff:P} yields
\begin{equation}\label{form:DeltaP}
    \begin{aligned}
        (\Delta P_{\mu,\eps}) U
    &=\sum_{j+k=-1}R\bracketLd{j}\Delta\L_{\mu,\eps}R\bracketLd{k}U-\sum_{j+k+l=-1}R\bracketLd{j}\Delta\L_{\mu,\eps}R\bracketLd{k}\Delta\L_{\mu,\eps} R\bracketLd{l}U+\errorTol
    \end{aligned}
\end{equation}
Next, for any vector $U\in\cV_{0,0}$, since the projection $P_{0,0}U=U$, we have
\begin{align*}
    P_{\mu,\eps} P_{0, 0}U + (1-P_{\mu,\eps})(1-P_{0, 0})U&=P_{\mu,\eps}U=(P_{0,0}+\Delta P_{\mu,\eps})U=(1+\Delta P_{\mu,\eps})U.
\end{align*}
Combining this with $\Delta\L_{\mu,\eps}=O(\eps,\mu) $, we find 
\begin{equation}
\label{expansion:cK}
    \begin{aligned}
    \cK_{\mu,\eps}U&=(1-(\Delta P_{\mu,\eps})^2)^{-1/2}(1+\Delta P_{\mu,\eps})U\\
    &=(1+\Delta P_{\mu,\eps}+\frac{1}{2}(\Delta P_{\mu,\eps})^2)U+\errorTol.
\end{aligned}
\end{equation}
To compute  $U_j^{\mu,\eps}=\cK_{\mu,\eps}U_j$, we use the coefficients $R_{[j]}$ (in \eqref{invL00:4},  \eqref{form:DeltaP}, and \eqref{expansion:cK}). After that, we compute the inner product $(\B_{\mu,\eps}U_j^{\mu,\eps},U_k^{\mu,\eps})$ to obtain entries of $\Lmat_{\mu,\eps}$. These heavy computations are carried out in the Mathematica notebook \verb|2 Lmat_expansion.nb|. We record the result in the following proposition.
\begin{prop}\label{prop:expandL} The coefficients  $\Lmat^{j,k}$ in the expansion \eqref{expansion:Lmat_formal} are given as follows.

$\eps^0$:
\begin{align}
\label{eq:lmat_coeff_00}
    \Lmat^{0,0}&=
\begin{bmatrix}
0 & 0 & 0 & 0 \\
 0 & 0 & 0 & 0 \\
 0 & 0 & 0 & \frac{-h c_0^4-4 \th^2 c_0^2+2 h \th c_0^2-h \th^2}{4 \th^2 c_0^2} \\
 0 & 0 & 0 & 0 
\end{bmatrix},\\
\Lmat^{1,0}&=
\begin{bmatrix}
\frac{i \left(h \th^2+\th-h\right) c_0^3}{\th \left(c_0^2+\th\right)} & 0 & 0 & 0 \\
 0 & \frac{i \left(h \th^2+\th-h\right) c_0^3}{\th \left(c_0^2+\th\right)} & 0 & 0 \\
 0 & 0 & \frac{i \left(-h c_0^2+2 \th c_0^2+h \th\right)}{2 \th c_0} & 0 \\
 0 & 0 & 0 & \frac{i \left(-h c_0^2+2 \th c_0^2+h \th\right)}{2 \th c_0} 
\end{bmatrix}.
\end{align}
\begin{equation}
\begin{aligned}
(\Lmat^{2,0})_{21}
&=-(\Lmat^{2,0})_{12}\\
&=\frac{c_0^3\left(
    \begin{aligned}
    &-h^2 c_0^4+h^2\th^2 c_0^4+2h^2\th^3 c_0^2-2h^2\th c_0^2\\
    &\quad-2h\th^3-h^2+h^2\th^2-\th^2+2h\th
    \end{aligned}
\right)}{(c_0^2+\th)^3},
\\[0.6em]
(\Lmat^{2,0})_{34}
&=\frac{h^3(c_0^2-\th)^2}{12\th^2 c_0^2},
\qquad
(\Lmat^{2,0})_{43}=h.
\end{aligned}
\end{equation}

$\eps^1$:
\begin{align}
    \Lmat^{0,1}&=
    \begin{bmatrix}
0 & 0 & 0 & \frac{2 \th^2 c_0^4-c_0^4-\th^2}{4 c_0^2} \\
 0 & 0 & 0 & 0 \\
 0 & \frac{2 \th^2 c_0^4-c_0^4-\th^2}{2 \th c_0 \left(c_0^2+\th\right)} & 0 & 0 \\
 0 & 0 & 0 & 0 
\end{bmatrix}.
\end{align}
\begin{equation}
\begin{aligned}
(\Lmat^{1,1})_{13}
&=-\frac{i\th(c_0^2+\th)}{2c_0},
\\[0.6em]
(\Lmat^{1,1})_{24}
&=\frac{i\left(
    \begin{aligned}
    &2\th^2 c_0^8-c_0^8-4h\th^4 c_0^6+2h\th^2 c_0^6+2h c_0^6\\
    &\quad-4\th c_0^6+2\th^4 c_0^4-8h\th^3 c_0^4-10\th^2 c_0^4\\
    &\quad+8h\th c_0^4-2h\th^4 c_0^2-4\th^3 c_0^2+2h\th^2 c_0^2-\th^4
    \end{aligned}
\right)}{8c_0^2(c_0^2+\th)^2},
\\[0.6em]
(\Lmat^{1,1})_{31}
&=\frac{i\left(
    \begin{aligned}
    &2\th^2 c_0^8-c_0^8-4h\th^4 c_0^6+2h\th^2 c_0^6+2h c_0^6\\
    &\quad-4\th c_0^6+2\th^4 c_0^4-8h\th^3 c_0^4-10\th^2 c_0^4\\
    &\quad+8h\th c_0^4-2h\th^4 c_0^2-4\th^3 c_0^2+2h\th^2 c_0^2-\th^4
    \end{aligned}
\right)}{4\th c_0(c_0^2+\th)^3},
\\[0.6em]
(\Lmat^{1,1})_{42}
&=-i.
\end{aligned}
\end{equation}

$\eps^2$:
\begin{equation}
\label{eq:lmat_coeff_02}
\begin{aligned}
(\Lmat^{0,2})_{12}
&=\frac{
    \begin{aligned}
    &-7\th^4 c_0^8+5\th^2 c_0^8-c_0^8-2\th^5 c_0^6\\
    &\quad+4\th^3 c_0^6-2\th c_0^6+3\th^4 c_0^4-3\th^2 c_0^4\\
    &\quad-2\th^5 c_0^2-2\th^3 c_0^2-\th^4
    \end{aligned}
}{4\th^2 c_0^5(c_0^2+\th)},
\\[0.6em]
(\Lmat^{0,2})_{34}
&=\frac{
    \begin{aligned}
    &-\th^2 c_0^{10}-c_0^{10}+8\th^5 c_0^8-8\th^3 c_0^8+6\th c_0^8\\
    &\quad-6\th^4 c_0^6+10\th^2 c_0^6-24\th^5 c_0^4+12\th^3 c_0^4\\
    &\quad-9\th^6 c_0^2+7\th^4 c_0^2+6\th^5
    \end{aligned}
}{16\th^2 c_0^2(c_0^2+\th)^2}.
\end{aligned}
\end{equation}
All other entries of $\Lmat^{2,0}$, $\Lmat^{1,1}$, and $\Lmat^{0,2}$ vanish.
\end{prop}
\subsection{Exact vanishing }
\def\nzOddVec{\vet{even}{odd}}
\def\nzEvenVec{\vet{odd}{even}}
\def\bracketMu#1{^{[#1]}}
 In this subsection, we prove that certain entries of the $\mu^0$ and $\mu^1$ coefficients of $\Lmat_{\mu, \eps}$ vanish exactly, i.e., without any finite order approximation in $\eps$. This will be obtained by leveraging parity of the (generalized) eigenvectors $U_{k, \sigma}$ of $\cL_{0, 0}$ and parity-preserving/switching property of various operators. To track parity, 
we denote by  $\nzOddVec$ the class of vectorial functions on $\T$ whose first component is even and second component is odd. The class  $\nzEvenVec$ is defined analogously. We say that an operator $T$ preserves parity if each of the above two classes is invariant under $T$. On the other hand, we say that $T$ switches parity if $T$ maps one class to the other.  Define $F\bracketMu{j}=F\bracketMu{j}(\eps)=[\mu^j]F$. It can be readily verified that  $\B\bracketMu{0}$ and  $\cK\bracketMu{0}$ preserve parity,  while $\B\bracketMu{1}$, $\L\bracketMu{0}$, and $J$ switch parity. A key observation is the following about $\cV_{0, \eps}$.
\begin{prop}\label{prop:evL0eps}
For all $\gamma \in \Rr$, there exists $\eps_0=\eps_0(\gamma)>0$ such that for all $0<|\eps|<\eps_0$, all  $ \nzEvenVec$ vectors in  the generalized eigenspace $\cV_{0, \eps}$ are eigenvectors of the eigenvalue $0$ of $\L_{0, \eps}$.
\end{prop}
\begin{proof}
We recall  Theorem \ref{thm:eigenvectors_R} for  the four (generalized) eigenvectors $R_j$ of eigenvalue $0$ of $\Ls_\eps$. We have $\Pi_{01}R_1=\Pi_{01}R_2=0$ since $(R_1)_1=0$ and $(R_2)_1=\p_x\eta^*$. Moreover,  \eqref{eq:solution_with_P} and \eqref{expand:Stokes_without_g_all_depth} imply
\[
(R_4)_1=\p_P\breve\eta\vert_{(a, P)=(\eps, 0)}=1+\p_h[\eta^*(x; \eps, h)]=1+O(\eps)>0
\]
for $\eps$ sufficiently small. Hence, 
 \[
 \Pi_{01}\big(R_3-CR_4\big)=0, \quad C:=(\Pi_{01}R_3)/(\Pi_{01}R_4).
 \]
Then, by virtue of Lemma \ref{prop:convert_RU}, the vectors 
 \[
U_1^\sharp=B_0R_1=R_1=\vet{0}{1},\quad U_2^\sharp=B_0 R_2,\quad U_3^\sharp=B_0(R_3-CR_4)
\]
are well-defined and satisfy 
\bq\label{relationRjUj}
 \L_{0,\eps} U^\sharp_j=B_0 \Ls_\eps  R_j-\vet{0}{s_j},\quad 1\le j\le 3,
\eq
where $s_j\in \Cc$ and $B_0$ is given by \eqref{def:B0}. From the formula \eqref{form:L0eps} of $\L_{0, \eps}$, we have  $\L_{0,\eps} U^\sharp_1=0$. Since $R_2=\nzEvenVec$, $B_0$ preserves parity, and $\L_{0, \eps}$ switches parity, we obtain  $U^\sharp_2=\nzEvenVec$ and $ \L_{0,\eps} U^\sharp_2=\nzOddVec$.  But $\Ls_\eps R_2=0$, so \eqref{relationRjUj} implies that $s_2=0$ and $\L_{0,\eps} U^\sharp_2=0$. Thus, $U^\sharp_1$ and $U^\sharp_2$ are eigenvectors of $\L_{0, \eps}$. 
 Next, we rewrite \eqref{relationRjUj} as
\[
 \L_{0,\eps} U^\sharp_j=B_0 \Ls_\eps  R_j-s_jU^\sharp_1.
\]
Recalling from \eqref{LR3R4} that $\Ls_\eps R_3$, $\Ls_\eps R_4\in \Span\{R_1, R_2\}$, we deduce 
\[
 \L_{0,\eps} U^\sharp_3\in \Span\{B_0R_1, B_0R_2, U^\sharp_1\}=\Span\{U^\sharp_1, U^\sharp_2\}.
 \]
Since $\vet{1}{0}=U_{0, +}$ is a generalized eigenvector of $\L_{0,0}$ with eigenvalue $0$, we have $U_4^\sharp:=\cK_{0,\eps} \vet{1}{0}\in\cV_{0,\eps}$. By Taylor expansion in $\eps\ne 0$, we find that  $U_1^\sharp,\dots,U_4^\sharp$ are linearly independent and hence form a basis of $\cV_{0, \eps}$. In this basis, $U^\sharp_1$, $U^\sharp_2\in \nzEvenVec$, and  $U^\sharp_3$, $U^\sharp_4\in \nzOddVec$ because  $R_3$, $R_4\in \nzOddVec$. Consequently, any  $ \nzEvenVec$ vector in $\cV_{0, \eps}$ must be spanned by $ U^\sharp_1$ and $U^\sharp_2$, and hence must be an eigenvector of eigenvalue $0$. 
\end{proof}
\begin{prop}
For all $\gamma \in \Rr$, there exists $\eps_0=\eps_0(\gamma)>0$ such that for all $|\eps|<\eps_0$, the following holds. There are analytic functions $b_j=b_j(\eps)\in\anaR{\eps}$ and $a_j=a_j(\mu,\eps)\in\anaR{\mu,\eps}$ such that
\begin{equation}
    \label{expansion:Lmat_formal_full}
    \begin{aligned}
        \Lmat_{\mu,\eps}
        &=\sum_{0\le j+k\le 2} \Lmat^{j,k}\mu^j\eps^k
        +\mu^3
\begin{bmatrix}
    \frac{i a_{11}}{\beta_1} & \frac{a_{12}}{\beta_1} & i a_{13} & a_{14} \\
    \frac{a_{15}}{\beta_1} & \frac{i a_{11}}{\beta_1} & a_{16} & i a_{17} \\
    \frac{i a_{17}}{\beta_1} & \frac{a_{14}}{\beta_1} & i a_{18} & a_{19} \\
    \frac{a_{16}}{\beta_1} & \frac{i a_{13}}{\beta_1} & a_{20} & i a_{18}
\end{bmatrix}+\mu ^2 \eps  \begin{bmatrix}
\frac{i a_1}{\beta _1} & \frac{a_2}{\beta _1} & i a_3 & a_4 \\
 \frac{a_5}{\beta _1} & \frac{i a_1}{\beta _1} & a_6 & i a_7 \\
 \frac{i a_7}{\beta _1} & \frac{a_4}{\beta _1} & i a_8 & a_9 \\
 \frac{a_6}{\beta _1} & \frac{i a_3}{\beta _1} & a_{10} & i a_8 
\end{bmatrix}\\
&+\mu  \eps ^2 \begin{bmatrix}
\frac{i b_{11}}{\beta _1} & \frac{b_{12}}{\beta _1} & i b_{13} & b_{14} \\
 0 & \frac{i b_{11}}{\beta _1} & 0 & i b_{17} \\
 \frac{i b_{17}}{\beta _1} & \frac{b_{14}}{\beta _1} & i b_{18} & b_{19} \\
 0 & \frac{i b_{13}}{\beta _1} & 0 & i b_{18} 
\end{bmatrix}+\eps ^3 \begin{bmatrix}
0 & \frac{b_2}{\beta _1} & 0 & b_4 \\
 0 & 0 & 0 & 0 \\
 0 & \frac{b_4}{\beta _1} & 0 & b_9 \\
 0 & 0 & 0 & 0 
\end{bmatrix},
    \end{aligned}
\end{equation}
where the coefficients $\Lmat^{j,k}$ are given in Proposition \ref{prop:expandL}. 
\end{prop}
\begin{proof}
When $\eps=0$, \eqref{expansion:Lmat_formal_full} follows from the analyticity in  $\mu$ of $\Lmat_{\mu, 0}$ and  Proposition \ref{prop:structure_Lmat}. We consider $\eps\ne 0$ in the remainder of this proof.  We recall that  $U\bracketMu{0}_{k,\pm}=\cK\bracketMu{0}U_{k,\pm}$, where $U_{k,+}=\nzOddVec$ and $U_{k,-}=\nzEvenVec$, $k=0, 1$. Consequently,  $U_{k,+}\bracketMu{0}=\nzOddVec$ and $U_{k,-}\bracketMu{0}=\nzEvenVec$. Then, Proposition \ref{prop:evL0eps} implies $\L\bracketMu{0}U_{k,-}\bracketMu{0}=0$, or equivalently, 
\bq\label{vanish:BUminus}
\B\bracketMu{0}U_{k,-}\bracketMu{0}=0.
\eq
Clearly $[\mu^0](\B_{\mu,\eps} U_{k,\sigma}^{\mu,\eps},U_{j,\tau}^{\mu,\eps})=(\B\bracketMu{0}U_{k,\sigma}\bracketMu{0},U_{j,\tau}\bracketMu{0})$. When at least one of $\sigma$ and $\tau$ is a minus sign, using \eqref{vanish:BUminus} and the fact that $\B\bracketMu{0}$ is a symmetric operator, we obtain $(\B\bracketMu{0}U_{k,\sigma}\bracketMu{0},U_{j,\tau}\bracketMu{0})=0$. Thus, from the formula \eqref{def:Lmat} of $\Lmat_{\mu, \eps}$, we conclude
\begin{equation}\label{Lmu0}
    [\mu^0]\Lmat=
    \begin{bmatrix}
        0 & * & 0 &*\\
        0 & 0& 0& 0\\
        0 & * & 0 &*\\
        0 & 0& 0& 0
    \end{bmatrix},
\end{equation}
where $``*"$ is simply a placeholder for the entry.

Next, we consider the coefficient 
\[
[\mu^1](\B_{\mu,\eps} U_{k,-}^{\mu,\eps},U_{j,-}^{\mu,\eps})=(\B\bracketMu{1}U_{k,-}\bracketMu{0},U_{j,-}\bracketMu{0})+(\B\bracketMu{0}U_{k,-}\bracketMu{1},U_{j,-}\bracketMu{0})+(\B\bracketMu{0}U_{k,-}\bracketMu{0},U_{j,-}\bracketMu{1})=:\sum_{\ell=1}^3I_\ell.
\]
$I_1=0$ because $\B\bracketMu{1}$ switches  parity, and $I_2=I_3=0$ thanks to \eqref{vanish:BUminus}. As a result, 
\begin{equation}\label{Lmu1}
    [\mu^1]\Lmat=
    \begin{bmatrix}
        * & * & * & *\\
        0 & * & 0 & *\\
        * & * & * & *\\
        0 & * & 0 & *
    \end{bmatrix}.
\end{equation}
From  \eqref{Lmu0}, \eqref{Lmu1}, and Proposition \ref{prop:structure_Lmat} for  $\Lmat_{\mu,\eps}$, we obtain the expansion \eqref{expansion:Lmat_formal_full}.
\end{proof}
\section{Spectral analysis}
In this section, we use the form \eqref{expansion:Lmat_formal_full} of the $4\times 4$ matrix  $\Lmat_{\mu, \eps}$ to find its spectrum.   We employ the direct method in our infinite-depth paper \cite{ChenNguyen}, which can be summarized as follows. After a possible rescaling of  the characteristic polynomial to factor out powers of parameters, this method begins with finding all roots in $i\Rr$  when $\mu=\eps=0$. The analytic implicit function theorem allows us to continue the simple roots as simple roots in $i\Rr$. A Newton-like iteration gives higher-order-in-$(\mu, \eps)$ approximations of these simple roots. After dividing by all linear factors given by simple roots, we end up with a lower-degree polynomial whose coefficients depend on the above simple roots. In the final step, we analyze roots of this reduced polynomial, which is a quadratic one in our case. 
\subsection{Technical note: Wolfram Mathematica}
The Mathematica codes in this section are selected from the companion notebook \verb|3 brute_force_eigenvalues.nb|. We implement a function that displays the coefficients of a polynomial \texttt{poly} in two variables \texttt{var1} and \texttt{var2}, in tabular form. If an entry contains more than $500$ characters, we replace it with $``\dots"$ in the output. This function was used in our earlier work \cite{ChenNguyen}. 
\begin{mmacode}
 \lstinline|showCoefficients[poly_, {var1_, var2_}, cutOff_:500]|
 
 \lstinline|showCoefficients[poly_]|
\end{mmacode}
 For example,
\begin{mmain}
    $\text{showCoefficients}(b \mu  \eps ^2+c \eps ^3+d \mu ^4+e+2^{10000} \mu ^6 \eps ^6+\mu ^2 \eps ,\{\eps ,\mu \})$
\end{mmain}
\begin{mmaout}
    \[
    \left(
\begin{array}{cccccccc}
 \text{} & 1 & \mu  & \mu ^2 & \mu ^3 & \mu ^4 & \mu ^5 & \mu ^6 \\
 1 & e & 0 & 0 & 0 & d & 0 & 0 \\
 \eps  & 0 & 0 & 1 & 0 & 0 & 0 & 0 \\
 \eps ^2 & 0 & b & 0 & 0 & 0 & 0 & 0 \\
 \eps ^3 & c & 0 & 0 & 0 & 0 & 0 & 0 \\
 \eps ^4 & 0 & 0 & 0 & 0 & 0 & 0 & 0 \\
 \eps ^5 & 0 & 0 & 0 & 0 & 0 & 0 & 0 \\
 \eps ^6 & 0 & 0 & 0 & 0 & 0 & 0 & \text{...} \\
\end{array}
\right)
    \]
\end{mmaout}
All the variables $b,c,d,e$ are put into the coefficients in corresponding positions, and the long coefficient $2^{10000}$ of $\eps^6\mu^6$ is replaced by $``\dots"$.

The analytic functions $a_{k}(\mu,\eps)$ and $b_k(\eps)$ in \eqref{expansion:Lmat_formal_full} are power series in $\eps$ and possibly $\mu$. However, the dependence on $\eps$ and $\mu$ is not explicitly coded in Mathematica, and the software treats them as independent variables. This causes no issue for our purpose, as our use of Mathematica is restricted to algebraic manipulations, every step of which is a polynomial in $a_k(\mu,\eps)$ and  $b_k(\eps)$. There is only one point to keep in mind when interpreting the output table: lower order terms with coefficients $a_{k}(\mu,\eps)$ and $b_k(\eps)$ may implicitly contribute to higher order terms. For example, consider the analytic function
\[
f_{example}(\mu,\eps)=\mu ^3 a_{1}(\mu,\eps)+\eps ^3 a_{2}(\mu,\eps)+5 \mu ^3 \eps ^3+3 \mu ^2 \eps ^2+1,
\]
which is displayed as 
\begin{mmaout}
    \[
    \left(
\begin{array}{ccccc}
 \text{} & 1 & \mu  & \mu ^2 & \mu ^3 \\
 1 & 1 & 0 & 0 & a_{1} \\
 \eps  & 0 & 0 & 0 & 0 \\
 \eps ^2 & 0 & 0 & 3 & 0 \\
 \eps ^3 & a_{2} & 0 & 0 & 5 \\
\end{array}
\right)
    \]
\end{mmaout}
From this table we may read that 
\[
[\mu^0\eps^0]f_{example}=1,[\mu^2\eps^2]f_{example}=3.
\]
On the other hand, the coefficient $[\mu^3\eps^3]f_{example}$ is not necessarily equal to $5$, since it may also be affected by $[\eps^3]a_1$ and $[\mu^3]a_2$.

\def\poly{\mathtt{p}}
\def\polyNew{\mathtt{p}_\mathtt{n}}
\def\phparen#1{\phi^{(#1)}}
\def\deri{\mathtt{d}}

\subsection{Characteristic polynomial }\label{subsect:char}
The eigenvalues of $\L_{\mu,0}$ are given by $i\omega_{k+\mu}^\pm,k\in\Z$ as in the formula \eqref{eq:omega_k_pm}. As $\mu\to 0$, the four eigenvalues near $0$ are $i\omega_{k+\mu}^\pm$ given by the following parameters. 
\begin{align*}
    \omega^\pm_\mu&= c_0\mu-\frac{\gamma}{2}h\mu \pm\mu\sqrt{h(1+\frac{\gamma^2}{4}h)}+O(\mu^2),\\
    \omega_{1+\mu}^-&=-\mu\frac{c_0^3}{\th}\frac{h-\th-h\th^2}{c_0^2+\th}+O(\mu^2),\\
    \omega_{-1+\mu}^+&=-\mu\frac{c_0^3}{\th}\frac{h-\th-h\th^2}{c_0^2+\th}+O(\mu^2).
\end{align*}
Consequently,  $\omega_{1+\mu}^--\omega_{-1+\mu}^+=O(\mu^2)$. To handle the square root in the above expansion of $\om_\mu^\pm$, we introduce the variable $H=\sqrt{\frac{1+\frac{\gamma^2}{4}h}{h}}$. In the symbolic calculation below, we replace every explicit occurrence of $h$ by $(H^2-\gamma^2/4)^{-1}$, except the $h$ in $\th$. Replacing $H$ by $-H$ interchanges the two roots $\om_\mu^\pm$. 

The characteristic polynomial of $\Lmat_{\mu,\eps}$ is
\begin{equation}
    \label{def:poly}
    \poly(\lambda;\mu,\eps)=\det(\lambda I-\Lmat_{\mu,\eps}),
\end{equation}
which is a quartic polynomial with coefficients in $\anaC{\mu,\eps}$. We make the change of variables
\begin{equation}
    \lambda=\mu f,
\end{equation}
and define a new polynomial in $f$,
\begin{equation}
    \label{def:polyNew}
    \polyNew(f;\mu,\eps):=\frac{\poly(\mu f;\mu,\eps)}{\mu^4}.
\end{equation}
We can check using Mathematica that $\polyNew$ contains no negative powers of $\mu$, so that 
\begin{equation}
     \polyNew(f;\mu,\eps)\in\anaC{\mu,\eps}[f].
\end{equation}
Moreover, we have
\begin{equation}
\label{eq:pn_coeff_real}
    [f^4]\polyNew=1,\quad [f^3]\polyNew,[f^1]\polyNew\in i\anaR{\mu,\eps},\quad [f^2]\polyNew,[f^0]\polyNew\in\anaR{\mu,\eps}.
\end{equation}
\subsection{First two roots of $\polyNew$}
To begin with, we let $\mu=\eps=0$ in $\polyNew$ and factorize the obtained polynomial $\polyNew(f;0,0)$ in Mathematica. The factorization yields two simple roots 
\[
\phi_\pm ^{(0)}:=i(c_0-\frac{\gamma}{2}h\pm h H),
\]
and a double root 
\[
\wt\phi=-i\frac{c_0^3}{\th}\frac{h-\th-h\th^2}{c_0^2+\th}.
\]
Clearly, $\phparen{0}_+\ne \phparen{0}_-$. Let us check that neither of them is equal to $\wt \phi$. If $\phparen{0}_\pm=\wt \phi$, then $\pm h H=-c_0+\frac{\gamma}{2}h-\frac{c_0^3}{\th}\frac{h-\th-h\th^2}{c_0^2+\th}$. Squaring both sides gives $h^2\frac{1+\frac{\gamma^2}{4}h}{h}=\Big(-c_0+\frac{\gamma}{2}h-\frac{c_0^3}{\th}\frac{h-\th-h\th^2}{c_0^2+\th}\Big)^2$.  But we have 
\begin{align*}
&h^2\frac{1+\frac{\gamma^2}{4}h}{h}-\Big(-c_0+\frac{\gamma}{2}h-\frac{c_0^3}{\th}\frac{h-\th-h\th^2}{c_0^2+\th}\Big)^2\\
&\qquad=-\frac{c_0^2 \left(c_0^2 h \th-hc_0^2 -h \th+h-\th\right) \left(c_0^2 h \th+c_0^2 h+h \th+h-\th\right)}{\left(c_0^2+\th\right){}^2}.
\end{align*}
 The two factors satisfy
\begin{equation}
    \label{eq:disc-denominator-signs}
    \begin{aligned}
   & c_0^2h\th-hc_0^2-h\th+h-\th
        =c_0^2h(\th-1)+(h-h\th-\th)<0,\\
    &c_0^2h\th+c_0^2h+h\th+(h-\th)
        >0.
    \end{aligned}
\end{equation}
Here we used $0<\th<1$ and
\[
    h-h\th-\th=-\frac{e^{2h}-2h-1}{e^{2h}+1}<0,
    \qquad h-\th>0.
\]
Thus the difference is strictly positive, which rules out $\phi_\pm^{(0)}=\wt\phi$. It follows that
\[
\deri:=\p_f \polyNew(\phparen{0}_+;0,0)=(\phparen{0}_+-\phparen{0}_-)(\phparen{0}_+-\wt\phi)^2\ne 0.
\]
Hence, the analytic implicit function theorem gives a unique analytic root
$\phi^*_+(\mu,\eps)$ near $\phparen{0}_+$ with
$\phi^*_+(0,0)=\phparen{0}_+$. For $(\mu, \eps)$ near $(0, 0)$, the
coefficient properties \eqref{eq:pn_coeff_real} imply 
\[
    \polyNew(-\overline{f};\mu,\eps)
    =\overline{\polyNew(f;\mu,\eps)}.
\]
Since $\phparen{0}_+\in i\Rr$, $\phi^*_+(0,0)=-\overline{\phi^*_+(0,0)}$,  uniqueness implies that
$\phi^*_+=-\overline{\phi^*_+}$ for all $(\mu,\eps)$ near origin, so $\phi^*_+\in i\Rr$.

By the variant of Newton's iteration described in \cite[Eq. (5.9)]{ChenNguyen}, we obtain a second-order approximation $\phparen{2}_+$ of $\phi^*_+$. The expression of $\phparen{2}_+$  is calculated and stored in the Mathematica notebook \verb|3 brute_force_eigenvalues.nb|. It consists of the $\mu^m \eps^n$ terms for $m+n\le 2$ with  pure imaginary coefficients. Invoking \cite[Lemma 5.2]{ChenNguyen},  we conclude that  $\phi^*_+(\mu,\eps)=\phparen{2}_+(\mu,\eps)+i r(\mu^3,\mu^2\eps,\mu\eps^2,\eps^3)$.

Since $\polyNew$ is invariant under
$H\mapsto-H$, this substitution interchanges the two simple roots. The root $\phi_-^*$ is obtained
from $\phi_+^*$ by replacing $H$ with $-H$. The product $\phi^*_+\phi^*_-$ and the sum  $\phi^*_+ + \phi^*_-$ are even in $H$, with exact expressions of terms of degree two or less, and with remainder terms of  $r(\mu^3,\mu^2\eps,\mu\eps^2,\eps^3)$ and $ir(\mu^3,\mu^2\eps,\mu\eps^2,\eps^3)$. We therefore  substitute $H^2=(1+\gamma^2h/4)/h$ to eliminate the variable $H$.

\subsection{The remaining two roots}
We denote  the remaining two roots of $\polyNew(f;\mu,\eps)$ by $\wt \phi_\pm(\mu,\eps)$, and let the quadratic equation that $\wt \phi_\pm$ satisfies be 
\begin{equation}
    \label{eq:quad_phi1_phi2}
    f^2+pf+q=0.
\end{equation}
Since
\begin{align*}
\phi_+^*+\phi_-^*+\wt\phi_++\wt\phi_-&=-[f^3]\polyNew,\\
\phi_+^*\phi_-^*+\phi_+^*\wt\phi_++\phi_+^*\wt\phi_-+\phi_-^*\wt\phi_++\phi_-^*\wt\phi_-+\wt\phi_+\wt\phi_-&=[f^2]\polyNew,
\end{align*}
we get
\[
p=-(\wt\phi_++\wt\phi_-)=[f^3]\polyNew +(\phi^*_++\phi^*_-)\in i\Rr
\]
and 
\begin{align*}
    q=\wt\phi_+\wt\phi_- &=[f^2]\polyNew-(\phi^*_++\phi^*_-)( \wt\phi_++\wt\phi_-)-\phi_+^*\phi_-^*\\
    &=[f^2]\polyNew+(\phi^*_++\phi^*_-)p-\phi_+^*\phi_-^*\in \Rr.
\end{align*}
The discriminant of the quadratic equation is $\Delta=p^2-4q\in \Rr$, and the two roots are
\bq
\wt\phi_\pm=\frac{-p\pm\sqrt\Delta}{2}.
\eq
We expand $p$ and $\Delta$ using Mathematica: 

\begin{mmain}
    $p\mathbin{//}\mathrm{showCoefficients}$
\end{mmain}
\begin{mmaout}
\[
\left(
\begin{array}{cccc}
    &1&\mu&\mu^2\\[0.3em]
1&
-\dfrac{2i(-c_0^3h+c_0^3\th+c_0^3h\th^2)}
       {\th(c_0^2+\th)}
&0&-\dfrac{4ic_0a_{11}}{\th(c_0^2+\th)}\\[0.6em]
\eps&0&\cdots&0\\
\eps^2&\cdots&0&0
\end{array}
\right)
\]
\end{mmaout}

\begin{mmain}
    $\mathit{quadraticDelta}\mathbin{//}\mathrm{showCoefficients}$
\end{mmain}
\begin{mmaout}
\[
\left(
\begin{array}{cccc}
    &1&\mu&\mu^2\\
1&0&0&\cdots\\
\eps&0&0&0\\
\eps^2&\cdots&0&0
\end{array}
\right)
\]
\end{mmaout}
The outputs read
\begin{align*}
    p
    &=-\frac{2i(c_0^3h\th^2-c_0^3h+c_0^3\th)}
             {\th(c_0^2+\th)}
      +ir(\mu^2,\mu\eps,\eps^2),\\
    \Delta
    &=\coeffM\mu^2+\coeffE\eps^2
      +r(\mu^3,\mu^2\eps,\mu\eps^2,\eps^3),
\end{align*}
where
\bq\label{def:cmu}
    \coeffM
    =-\frac{4c_0^6}{(c_0^2+\th)^6}
      \Big[\th^4+(1-\th^2)
      \big(h^2c_0^4+2h^2\th c_0^2+(h-\th)^2\big)\Big]^2,
\eq
\bq\label{def:ceps}
\begin{aligned}
    \coeffE
    &=\frac{\Omega}{c_0^2\th^2(c_0^2+\th)^4}
    \frac{\th^4+(1-\th^2)
          \big(h^2c_0^4+2h^2\th c_0^2+(h-\th)^2\big)}
         {\big[c_0^2h\th-hc_0^2-h\th+h-\th\big]
          \big[c_0^2h\th+c_0^2h+h\th+(h-\th)\big]}.
\end{aligned}
\eq
Here $\Omega=\Os(h, d)$,  $d:=c_0^2$, is the following sextic polynomial in $d$:
\begin{equation}
    \label{eq:Omega_explicit}
    \Os(h,d)=\sum_{j=0}^{6}\Oscoeff{j}(h)d^j,
\end{equation}
\bq\label{form:Omegaj}
\begin{aligned}
\Oscoeff{0}(h)
    &=\th^4\big[(1-\th^2)(h-\th)^2+\th^4\big],\\
\Oscoeff{1}(h)
    &=2\th^3\big[(h-\th)^2-h\th^3+3\th^4-h^2\th^4\big],\\
\Oscoeff{2}(h)
    &=\th^2\big[3(h-\th)^2+h^2\th^2(4\th^2-7)
              +2h\th^3(3-5\th^2)+9\th^4\big],\\
\Oscoeff{3}(h)
    &=2\th\big[(h-\th)^2+h^2\th^2(3\th^2-4)
              +h\th^3(3-5\th^2)+\th^4(6-\th^2)\big],\\
\Oscoeff{4}(h)
    &=(h-\th)^2+h^2\th^2(-9+18\th^2-10\th^4)\\
    &\quad{}+2h\th^3(5-9\th^2+\th^4)+\th^4(7-\th^2),\\
\Oscoeff{5}(h)
    &=2\th\big[-h^2(1-\th^2)^3+h\th(1-\th^2-\th^4)
              +2\th^2(1-\th^2)\big],\\
\Oscoeff{6}(h)
    &=h\big[2\th(1-3\th^2+3\th^4)
             -h(1-\th^2)(1-5\th^2+7\th^4)\big].
\end{aligned}
\eq
The common factor $\th^4+(1-\th^2)
          \big(h^2c_0^4+2h^2\th c_0^2+(h-\th)^2\big)$ in the numerators of  $\coeffM$ and $\coeffE$ is positive since $0<\th<1$. By
\eqref{eq:disc-denominator-signs}, the two bracketed factors in the
denominator of $\coeffE$ have opposite signs. We deduce
\bq
\coeffM<0,\quad \sign(\coeffE)=-\sign(\Omega).
\eq
Combining the considerations in this and the preceding subsections, we have proven \ref{theo:main1}, where we note that the two purely imaginary roots $ \ld^*_\pm$ are given by
\bq
 \ld^*_\pm=\mu \phi^{(2)}_\pm+i\mu r_\pm^*(\mu^3, \mu^2\eps, \mu\eps^2, \eps^3).
 \eq
\subsection{Analysis of $\Om$}
In this subsection, we use  formal proofs verified in Rocq to establish Theorem \ref{theo:main2}, which gives a global sign classification of  $\Omega$. For the reader's convenience, we restate Theorem \ref{theo:main2} below.
\begin{theo}[The sign of $\Omega$]
\label{thm:Omega_sign}
The following statements hold.
\begin{enumerate}
    \item
    There is a unique point
    $(\hcrit,\dcrit)\in(0,\infty)^2$
    satisfying $\Omega_\sharp(\hcrit,\dcrit)
   =\partial_d\Omega_\sharp(\hcrit,\dcrit)=0$. 
    Moreover,
    \[
        1.36271<\hcrit<1.36272,
        \qquad
        0.86010<\dcrit<0.86011.
    \]
  Denote 
    \[
        \gammacrit
        :=\frac{\sqrt{\dcrit}}{\tanh(\hcrit)}
          -\frac{1}{\sqrt{\dcrit}}
        \in(-0.02081,-0.02079).
    \]
     \item
    If $0<h<\hcrit$, then $ \Omega(h,\gamma)>0$ for all $\gamma\in\mathbb R$. At the critical depth, $\Omega(\hcrit,\gamma)\ge 0$ for all $\gamma\in\Rr$, with equality if and only if $\gamma=\gammacrit$.
    \item Consider $h>\hcrit$. The equation
    $\Omega(h, \cdot)=0$ has exactly two real roots,
    denoted by $\gamma_l(h)<\gamma_r(h)$, which are smooth functions of $h$. The sign of $\Omega$ is given by
    \[
    \begin{cases}
        \Omega(h,\gamma)>0
        & \text{if }\gamma<\gamma_l(h)
          \text{ or }\gamma>\gamma_r(h),\\
        \Omega(h,\gamma)=0
        & \text{if }\gamma=\gamma_l(h)
          \text{ or }\gamma=\gamma_r(h),\\
        \Omega(h,\gamma)<0
        & \text{if }\gamma_l(h)<\gamma<\gamma_r(h).
    \end{cases}
    \]
    \item For any $\gamma\ge \frac{2}{\sqrt3}$ and any finite depth $h$, $\Om(h,\gamma)>0$. 
\end{enumerate}
\end{theo}
By virtue of (4), any constant vorticity  $\gamma\ge\frac{2}{\sqrt3}$ suppresses the Benjamin--Feir instability at any finite depth. For a depth-dependent lower bound, the Rocq proof script gives 
\bq
\gamma\ge \sqrt3-\frac{1}{\sqrt 3\tanh(h)}.
\eq
  Figure~\ref{fig:omega-contour} shows the zero set of $\Omega$ in the
$(h,\gamma)$-plane. The lower branch $\gamma_l(h)$ and the upper
branch $\gamma_r(h)$ meet at $A=(\hcrit,\gammacrit)$. 
\begin{rema} Based on the graph, we make the following observations, which  have not been rigorously proven. In the plot, the lower branch decreases with $h$, while the upper
branch has a local maximum at $B$ and a local minimum at $C$.
The upper branch increases up to
$B$, decreases from $B$ to $C$, and increases again beyond $C$.
For a vorticity strictly between the two horizontal levels marked
in the second panel, the line of constant vorticity intersects the
zero set of $\Omega$ at three depths $h_1<h_2<h_3$. Along this line,
$\Omega<0$ for $h\in(h_1,h_2)\cup(h_3,\infty)$, so the
Benjamin--Feir instability occurs when the depth is in these two intervals.
At $h=h_j$, $j=1,2,3$, our second-order criterion based on $\sign(\Omega)$ is inconclusive.
\begin{figure}[H]
    \centering
    \includegraphics[width=0.85\linewidth]{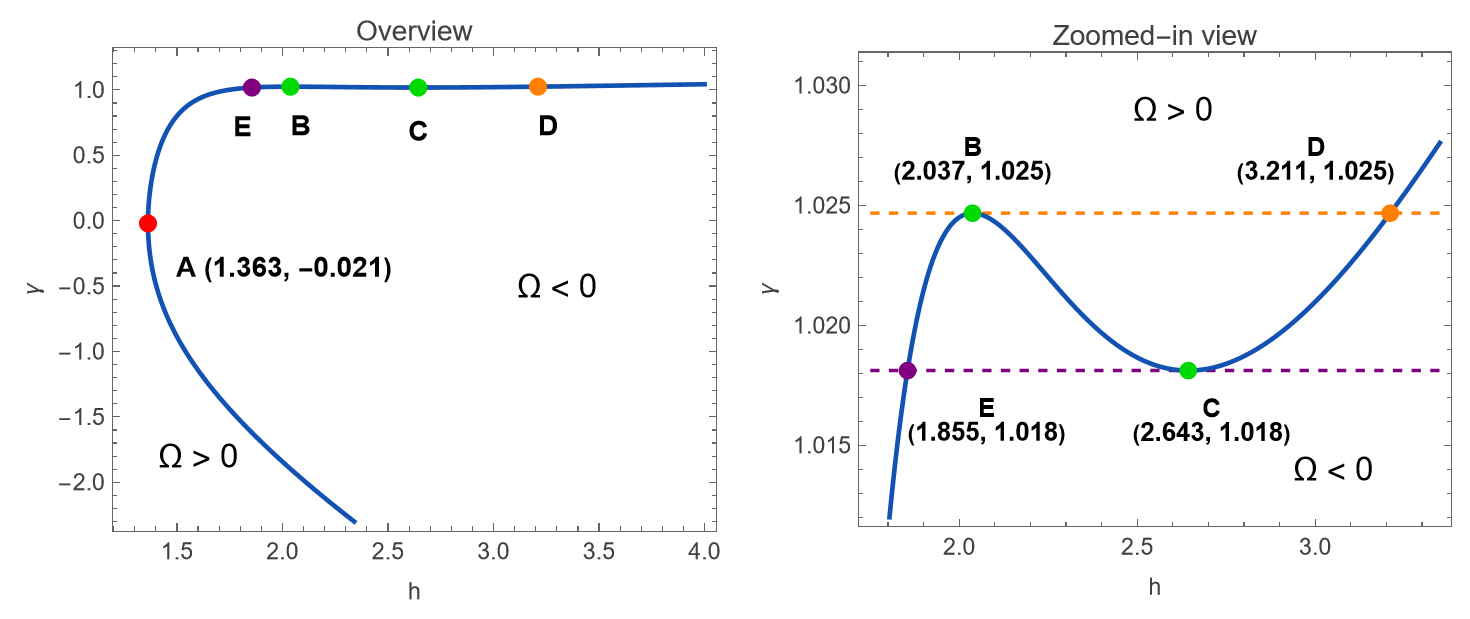}
    \caption{The zero set of $\Omega$ in the $(h,\gamma)$-plane.}
    \label{fig:omega-contour}
\end{figure}\label{figure:Om}
\end{rema}

\begin{proof}[Proof of theorem  \ref{thm:Omega_sign}]
    
For fixed $h>0$, the map $d\mapsto \sqrt d/\th-1/\sqrt d$ is
a strictly increasing bijection from $(0,\infty)$ onto $\mathbb R$.
Thus it suffices to study  positive roots and the sign of
$\Os(h,\cdot)$ in \eqref{eq:Omega_explicit}.

We record the main theorem statement from the Rocq proof script \verb|OmegaRootBifurcation.v|. Here \texttt{Omega} denotes
$\Os(h,d)$ from \eqref{eq:Omega_explicit}, and \texttt{Omega\_d} denotes
its derivative with respect to $d$. The root counts are for distinct
positive roots.
\begin{lstlisting}[
    basicstyle=\small\ttfamily,
    columns=fullflexible,
    keepspaces=true,
    breaklines=true,
    label={lst:rocq-root-bifurcation}
]
Theorem certified_direct_positive_root_bifurcation :
  (forall h d : R, is_derive (Omega h) d (Omega_d h d)) /\
  exists h0 d0,
    (0 < h0 /\ 0 < d0 /\
     Omega h0 d0 = 0 /\ Omega_d h0 d0 = 0) /\
    136271/100000 < h0 < 136272/100000 /\
    86010/100000 < d0 < 86011/100000 /\
    (forall h d,
      0 < h -> 0 < d ->
      Omega h d = 0 -> Omega_d h d = 0 ->
      h = h0 /\ d = d0) /\
    (forall h, 0 < h < h0 -> no_positive_root (Omega h)) /\
    exactly_one_positive_root (Omega h0) /\
    (forall h, h0 < h -> exactly_two_positive_roots (Omega h)).
\end{lstlisting}
Let us explain the code. In Rocq, \texttt{->} denotes implication and associates to the right, so
\texttt{A -> B -> C} means \texttt{A -> (B -> C)}, or equivalently, if propositions $A$ and $B$ both hold, then the proposition $C$ is true. \texttt{/\textbackslash} denotes logical conjunction $\land$. \texttt{Omega h d} applies \texttt{Omega} to the arguments \texttt{h} and \texttt{d}, and represents $\Omega_\sharp(h,d)$. It is proven that the critical point $(\hcrit,\dcrit)$ where $\Os(\hcrit,\dcrit)=\p_d\Os(\hcrit,\dcrit)=0$ with $\hcrit>0,\dcrit>0$ is unique. In addition, $\hcrit\in(1.36271,1.36272)$ and $\dcrit\in(0.86010,0.86011)$. It is also proven that positive root counts of the polynomial $\Os(h,\cdot)$ are $0,1,2$ for $h<\hcrit,h=\hcrit,h>\hcrit$ respectively. The uniqueness and the decimal bounds of $(\hcrit,\dcrit)$ are proven rigorously by some exact inequalities, without any floating point errors. These inequalities in Rocq can be considered coming from Taylor expansions in $(h,d)$, with explicit bounds of error terms. 

The smoothness of $\gamma_l(h)$ and $\gamma_r(h)$ comes from implicit function theorem, which is valid on the zero set $\Os(h,d)=0$ for $h>\hcrit$. The sign of $\Omega$ is based on the positivity of constant term in the sextic polynomial $\Os(h,\cdot)$ given in \eqref{form:Omegaj}.

Finally, in the Rocq proof script \verb|OmegaGammaPositivity.v|, it proves $\Om(h,\gamma)>0$ for any $\gamma\ge \sqrt3-\frac{1}{\sqrt 3\tanh(h)}$. The last assertion of the theorem follows from the fact that $\sqrt3-\frac{1}{\sqrt 3\tanh(h)}<\frac{2}{\sqrt3}$ for any $h>0$.

\begin{lstlisting}[
    basicstyle=\small\ttfamily,
    columns=fullflexible,
    keepspaces=true,
    breaklines=true
]
Definition gamma_star (h : R) : R :=
  sqrt 3 - 1 / (sqrt 3 * tanh h).

Theorem Omega_positive_from_gamma_star :
  forall h c0 gamma : R,
  0 < h -> 0 < c0 ->
  gamma = c0 / tanh h - 1 / c0 ->
  gamma_star h <= gamma ->
  0 < Omega h (c0^2).
\end{lstlisting}
\end{proof}

\appendix

\section{Expansion of the conformal mapping}
\label{appendix:expansion_conformal_mapping}
In this appendix, we use \eqref{eq:fixed_point_zeta_fe} to derive the expansions for $\zeta$ and $\fe$ in Lemma \ref{lemm:zeta}. We will invoke the following formula from \cite[formula (2.14)]{Craig_Sulem}:
\begin{equation}
\label{eq:expansion:G_order_eps_3}
\begin{aligned}
G(\eta) \psi 
& =\underbrace{G_0 \psi}_{O(\eps)}
+ \underbrace{\left( \eta D^2 - \eta_x \partial_x  \right) \psi -G_0(\eta G_0(\psi))}_{O(\eps^2)}\\
&+\underbrace{\left((\eta_x)^2G_0\psi+G_0\big(\eta G_0(\eta G_0\psi)\big)-\frac{1}{2}G_0(\eta^2D^2\psi)-\eta D^2(\eta G_0\psi)+\frac{1}{2}\eta^2 D^2G_0\psi\right)}_{O(\eps^3)}
+ \rem{4},
\end{aligned}
\end{equation}
where $G_0=|D|\tanh(h|D|)$ is the unperturbed Dirichlet to Neumann operator. 

When $\eps=0$, we have $\zeta(x)=x,\fe=0$, and hence $\fe=O(\eps)$. We write 
\begin{equation}
    \zeta(x)=x+\eps\zeta_1(x)+\eps^2\zeta_2(x)+\rem 3,\quad  \eta^*(x)=\eps\eta_1(x)+\eps^2\eta_2(x)+\rem 3,
\end{equation}
where $\eta_1=\th\cos(x)$ and $\eta_2$ are given in \eqref{expand:Stokes_without_g_all_depth}.

We first expand the equation of $\zeta(x)$ in \eqref{eq:fixed_point_zeta_fe} to first  order in $\eps$. Using that $\eta^*=O(\eps)$, we find
\begin{align*}
    x+\eps\zeta_1(x)+\rem 2&=x+\frac{\mathcal{H}}{\tanh\big(h|D|+O(\eps)\big)}\eta^*(x+O(\eps))=x+\frac{\mathcal{H}}{\tanh(h|D|)}\eta^*(x)+O(\eps^2).
\end{align*}
Consequently, 
\begin{equation}
\label{eq:zeta_1}
    \begin{aligned}
    \zeta_1(x)=\frac{\mathcal{H}}{\tanh(h|D|)}\eta_1(x)=\frac{\mathcal{H}}{\tanh(h|D|)}[\th\cos(x)]
    =\frac{\mathcal{H}}{\th}[\th\cos(x)]
    =\sin(x).
\end{aligned}
\end{equation}
On the other hand, the second equation in \eqref{eq:fixed_point_zeta_fe} gives
\begin{align*}
    \fe&=\fint\eta^*\Big(x+\eps\zeta_1(x)+O(\eps^2)\Big)dx=\fint\Big[\eta^*(x)+\eps\zeta_1(x)\p_x\eta^*(x)+\rem 3 \Big]dx.
\end{align*}
Inserting $\eta^*$ from \eqref{expand:Stokes_without_g_all_depth} and using \eqref{eq:zeta_1}, we deduce
\begin{align*}
    \fe&=\frac{\eps^2}{4}(c_0^2-2\th +c_0^2\th^2)+\fint \big(\eps\sin(x)\big)\big(-\eps\th\sin(x)\big)dx+\rem3\\
       &=\frac{\eps^2}{4}(c_0^2-4\th +c_0^2\th^2)+\rem 3.
\end{align*}
Next, we expand the first equation in \eqref{eq:fixed_point_zeta_fe} to second order in $\eps$:
\begin{align*}
    x+\eps\zeta_1(x)+\eps^2\zeta_2(x)+\rem 3&=x+\frac{\mathcal{H}}{\tanh\big(h|D|+O(\eps^2)\big)}\eta^*(x+\eps\sin(x)+O(\eps^2))\\
    &=x+\frac{\mathcal{H}}{\tanh\big(h|D|\big)}\Big[\eta^*(x)+\eps\zeta_1(x)\p_x\eta^*(x)\Big]+\rem 3.
\end{align*}
This yields 
\begin{align*}
    \zeta_2(x)&=\frac{\mathcal{H}}{\tanh\big(h|D|\big)}\Big[\eta_2(x)+\zeta_1(x)\p_x\eta_1(x) \Big]\\
    &=\frac{\mathcal{H}}{\tanh\big(h|D|\big)}\Big[C+\frac{1}{4}\Big(\frac{1}{\th}+\frac{c_0^2+\th(1-2c_0^4)+c_0^2\th ^2}{c_0^4}\Big)\cos(2x)-\th\sin^2(x)\Big]\\
    &=\frac{\mathcal{H}}{\tanh\big(h|D|\big)}\Big[\frac{1}{4}\Big(\frac{1}{\th}+\frac{c_0^2+\th(1-2c_0^4)+c_0^2\th ^2}{c_0^4}\Big)\cos(2x)+C-\frac{\th}{2}+\frac{\th}{2}\cos(2x)\Big]\\
    &=\frac{1}{4}\Big(\frac{1}{\th}+\frac{c_0^2+\th+c_0^2\th ^2}{c_0^4}\Big)\frac{\mathcal{H}}{\tanh\big(h|D|\big)}\cos(2x)\\
   &=\frac{1}{4}\Big(\frac{1}{\th}+\frac{c_0^2+\th+c_0^2\th ^2}{c_0^4}\Big)\frac{1}{\tanh(2h)}\sin(2x)\\
   &=\frac{1}{4}\Big(\frac{1}{\th}+\frac{c_0^2+\th+c_0^2\th ^2}{c_0^4}\Big)\frac{1+\th^2}{2\th}\sin(2x).
\end{align*}
\subsection*{Acknowledgements}  HQN was partially supported by NSF CAREER Grant DMS-2541807. The authors thank M. Berti for a detailed discussion on October 29, 2026, concerning the results of this paper and \cite{ChenNguyen}.

    \bibliographystyle{amsplain} 
    \bibliography{references}

\end{document}